\documentclass[a4paper,11pt]{article}

\usepackage[applemac]{inputenc}
\usepackage{enumerate}
\usepackage{lmodern}
\usepackage[T1]{fontenc}
\usepackage{verbatim}
\usepackage{textcomp}
\usepackage[english]{babel}
\usepackage[a4paper,vmargin={3.5cm,3.5cm},hmargin={2.5cm,2.5cm}]{geometry}
\usepackage[font=sf, labelfont={sf,bf}, margin=1cm]{caption}
\usepackage[pdftex]{hyperref}
\usepackage[pdftex]{color,graphicx}

\usepackage{amsmath,amsfonts,amssymb,amsthm,mathrsfs}

\usepackage{cleveref}
  \crefname{theorem}{Theorem}{Theorems}
  \crefname{lemma}{Lemma}{Lemmas}
  \crefname{remark}{Remark}{Remarks}
  \crefname{proposition}{Proposition}{Propositions}
  \crefname{definition}{Definition}{Definitions}
  \crefname{corollary}{Corollary}{Corollaries}
  \crefname{section}{Section}{Sections}
  \crefname{figure}{Figure}{Figures}

\usepackage{color}

\newtheorem{theorem}{Theorem}[]
\newtheorem{definition}[theorem]{Definition}
\newtheorem{proposition}[theorem]{Proposition}
\newtheorem{lemma}[theorem]{Lemma}

\theoremstyle{definition}

\def\llbracket{[\hspace{-.10em} [ }
\def\rrbracket{ ] \hspace{-.10em}]}

\def\noi{\noindent}

\def\t{{\mathcal T}}
\def\n{{\mathcal N}}

\def\ve{\varepsilon}
\def\wt{\widetilde}
\def\wh{\widehat}
\def\bm{{\bf m}}

\def\la{\longrightarrow}

\def\R{{\mathbb R}}
\def\P{{\mathbb P}}
\def\E{{\mathbb E}}
\def\N{{\mathbb N}}

\def\W{{\mathcal W}}
\def\z{\mathcal{Z}}
\def\w{\mathrm{w}}

\def\ov{\overline}
\def\build#1_#2^#3{\mathrel{
\mathop{\kern 0pt#1}\limits_{#2}^{#3}}}
\def\rem{\noindent{\bf Remark. }}

\title{The Brownian cactus II. \\ Upcrossings and local times of super-Brownian motion}
\author{Jean-Fran\c cois Le Gall\thanks{Universit\'e Paris-Sud}}
\date{}
\begin{document}

\maketitle

\begin{abstract}
We study properties of the random metric space called the Brownian map. For every 
$r>0$, we consider the connected components of the complement of the open ball of radius $r$ centered at the root,
and we let $\mathbf{N}_{r,\ve}$ be the number of those connected components that intersect the complement
of the ball of radius $r+\ve$. We then prove that $\ve^3\mathbf{N}_{r,\ve}$ converges as $\ve\to 0$
to a constant times the density at $r$ of the profile of distances from the root. In terms of
the Brownian cactus, this gives asymptotics for the number of vertices at height $r$ that have
descendants at height $r+\ve$. Our proofs are based on a similar approximation result 
for local times of super-Brownian motion by upcrossing numbers. Our arguments make
a heavy use of the Brownian snake and its special Markov property. 
\end{abstract}

\section{Introduction}

This paper is devoted to certain properties of the random metric space known as
the Brownian map, which can be viewed as a canonical model of random 
geometry in two dimensions. These properties are closely related to an approximation result
for local times of super-Brownian motion in terms of upcrossing numbers, which is
similar to the classical result for linear Brownian motion.

In order to present our main results, let $(\bm_\infty,D)$ denote the Brownian map. This is a 
random compact metric space, which is a.s.homeomorphic to the two-dimensional sphere and
has recently been shown to be the scaling limit in distribution, in the Gromov-Hausdorff sense, of several classes 
of random planar maps \cite{AA,BLG,LUniqueness,Mi-quad}. The Brownian map is equipped with a
volume measure $\lambda$, which in a sense is the uniform probability measure on $\bm_\infty$, and a distinguished point,
which we denote here by $\rho$. This point plays no particular role in the sense that,
if we ``re-root'' the Brownian map at another point $\tilde\rho$ chosen according to $\lambda$, the pointed metric spaces $(\bm_\infty, D, \rho)$
and $(\bm_\infty,D,\tilde\rho)$ have the same distribution \cite[Theorem 8.1]{LGeodesic}. For every $h>0$, let $B_h(\rho)$
stand for the open ball of radius $h$ centered at $\rho$. Then, on the event where $B_h(\rho)^c\not =
\varnothing$, $B_h(\rho)^c$ will have infinitely many connected components, but a compactness 
argument shows that only finitely 
many of them intersect $B_{h+\ve}(\rho)^c$, for any fixed $\ve>0$. Our first objective is to get 
precise information about the number of these components.  Recall that
the profile of distances from $\rho$ in $\bm_\infty$  is the 
probability measure $\Delta$  on $\R_+$ defined by
$$\Delta(A):=\int \lambda(\mathrm{d}x)\,\mathbf{1}_A(D(\rho,x)),$$
for any Borel subset $A$ of $\R_+$.
The measure $\Delta$ has a.s.~a continuous density with
respect to Lebesgue measure. 

 \begin{theorem}
 \label{connec-compo}
 For every $h>0$ and $\ve>0$, let $\mathbf{N}_{h,\ve}$ be the number of connected components of 
 $B_h(\rho)^c$ that intersect $B_{h+\ve}(\rho)^c$. Then,
 \begin{equation}
\label{intro-1}
 \ve^3\,\mathbf{N}_{h,\ve} \build{\la}_{\ve \to0}^{} \frac{c_1}{2}\,\mathbf{L}^h
 \end{equation}
 in probability. Here $\mathbf{L}^h$ is the
density at $h$ of the profile of distances from $\rho$ in $\bm_\infty$, and the constant $c_1>0$ is determined by the
identity  
$$\int_0^\infty \frac{\mathrm{d}u}{\sqrt{\frac{8}{3}u^3 + (c_1)^2}} = 1.$$
 \end{theorem}

Theorem \ref{connec-compo} can be reformulated in terms of the Brownian cactus
discussed in \cite{CLM}. Recall that, with any pointed geodesic compact metric space, one can
associate a rooted $\R$-tree called the cactus of the initial space. Roughly speaking, 
the root of the cactus corresponds to the distinguished point in the original space, and distances
from this point are in a sense preserved in the cactus. Furthermore, the
vertices of the cactus at a given height $h$, that is, at distance $h$ from the root, correspond to the connected components 
of the complement of the open ball of radius $h$ centered at the distinguished point (see
\cite[Section 2.5]{CLM}). The cactus associated with the Brownian map is called the
Brownian cactus (one of the main reasons for introducing this object is the fact that the convergence in distribution 
of discrete cactuses associated with random planar maps toward the Brownian cactus
has been proved in great generality \cite{CLM}). The quantity $\mathbf{N}_{h,\ve}$
is then equal to the number of vertices of the Brownian cactus at height $h$
that have descendants at height $h+\ve$, and \eqref{intro-1} shows that this number is
typically of order $\ve^{-3}$ when $\ve$ tends to $0$. 

Perhaps more surprisingly, the convergence \eqref{intro-1} is also closely related to an approximation result
for local times of super-Brownian motion in terms of upcrossing numbers. If $w:[0,T]\la \R$
is a continuous function defined on the interval $[0,T]$, and $h\in \R$, we say that $r\in[0,T)$ is an upcrossing time of 
$w$ from $h$ to $h+\ve$ if $w(r)=h$ and if there exists $t\in(r,T]$ such that $w(t)=h+\ve$
and $w(s)>h$ for every $s\in(r,t]$. Then, if $\mathrm{N}_{h,\ve}(T)$ is the number  
of upcrossing times from $h$ to $h+\ve$ of a standard linear Brownian motion $B$
over the time interval $[0,T]$, $(2\ve)^{-1}\,\mathrm{N}_{h,\ve}(T)$ converges a.s.~as $\ve\to 0$
to the local time of $B$ at level $h$ and at time $T$. This is the classical approximation of Brownian local times
by upcrossing numbers (see \cite[Section 2.4]{IM}
or \cite[Theorem VI.1.10]{RY}).   In view of a similar result for super-Brownian
motion, we would like to count upcrossing times for all ``historical paths'', and for this we need
introduce the historical super-Brownian motion (see \cite{DP,Dy} for the general theory
of historical superprocesses).

So let $\mathbf{Y}=(\mathbf{Y}_t)_{t\geq 0}$
be a one-dimensional historical super-Brownian motion. For every $t\geq 0$, $\mathbf{Y}_t$
is a random measure on the space $C([0,t],\R)$ of all continuous functions from $[0,t]$ into $\R$.
Informally, the support of $\mathbf{Y}_t$ consists of the historical paths followed between times $0$ and $t$ by
all ``particles'' alive at time $t$. The associated super-Brownian motion $\mathbf{X}=(\mathbf{X}_t)_{t\geq 0}$
is obtained from $\mathbf{Y}$ by the formula
$$\mathbf{X}_t(A)= \int \mathbf{Y}_t(\mathrm{d}w)\,\mathbf{1}_A(w(t)),$$
for any Borel subset $A$ of $\R$ and every $t\geq 0$. 
Then, for every $t>0$, 
$\mathbf{X}_t$ has a continuous density denoted by $u_t$ (see e.g.~\cite[Theorem III.4.2]{Pe}), and we
set, for every $x\in\R$,
$$\mathrm{L}^x=\int_0^\infty \mathrm{d}t\, u_t(x).$$
Clearly, the function $x\to \mathrm{L}^x$ is also the density of the occupation measure $\int_0^\infty \mathrm{d}t\,\mathbf{X}_t$, and for this reason
we call $\mathrm{L}^x$ the local time of $\mathbf{X}$ at level $x$. Note that these local times also exist in
dimensions $2$ and $3$ even though the measures $\mathbf{X}_t$ are then singular (see \cite{Fl,Su}).

We say that $r\geq 0$ is an upcrossing time of $\mathbf{Y}$ from $h$ to $h+\ve$ 
if there exist $t>r$ and a function $w\in C([0,t],\R)$ that belongs to the topological support of $\mathbf{Y}_t$,
such that $r$ is an upcrossing time of $w$ from $h$ to $h+\ve$.

\begin{theorem}
\label{SBM}
Assume that $\mathbf{X}_0=a\,\delta_0$
for some $a>0$, where $\delta_0$ denotes the Dirac measure at $0$. 
Let $h\in\R\backslash\{0\}$ and for every $\ve>0$, let $\mathscr{N}_{h,\ve}$ be the
number of upcrossing times of $\mathbf{Y}$ from $h$ to $h+\ve$. Then
\begin{equation}
\label{intro-2}
\ve^3\,\mathscr{N}_{h,\ve} \build{\la}_{\ve \to 0}^{} \frac{c_1}{2}\,\mathrm{L}^h,
\end{equation}
in probability. Here, $\mathrm{L}^h$ is the local time of $\mathbf{X}$ at level $h$, and the constant $c_1$ was defined in Theorem \ref{connec-compo}.
\end{theorem}

\rem  The definition of upcrossings in the superprocess setting can also be interpreted in terms
of the genealogical structure of super-Brownian motion. Recall that the
genealogy of $\mathbf{X}$ is coded by a random $\R$-tree, or more precisely by a countable collection
of random $\R$-trees. Each vertex of these trees is assigned a spatial location in $\R$, and the 
measure $\mathbf{X}_t$ is in a sense ``uniformly spread'' over the spatial locations of vertices at height $t$. Then upcrossing times of $\mathbf{Y}$ from $h$ to $h+\ve$ are in one-to-one correspondence with
vertices $v$ whose spatial location is equal to $h$ and which have (at least)
one descendant $v'$ with spatial location  $h+\ve$, such that spatial locations stay greater than $h$ on the line segment between $v$ and $v'$
in the tree. See Section \ref{sec-upcrossing} below for a rigorous presentation of this interpretation. 

\medskip
Our proof of Theorem \ref{connec-compo} relies on a version of the convergence of Theorem \ref{SBM}
under the excursion measure of super-Brownian motion (\cref{appro-LT}). Let us
explain the connection between connected components of the complement 
of a ball in the Brownian map and upcrossings of super-Brownian motion. We first recall that the
Brownian map is constructed as a quotient space of Aldous' Continuous Random Tree (the so-called CRT) for
an equivalence relation which is defined in terms of Brownian labels assigned to the
vertices of the CRT (see Section \ref{sec-applimap} for more details). Note that the CRT is just a conditional version of the random trees coding the genealogy of super-Brownian motion, and that the Brownian labels  
can be viewed as spatial locations in the superprocess setting. From the properties of the Brownian map, it is not too hard to prove that
connected components of $B_h(\rho)^c$ in $\bm_\infty$ correspond 
to connected components of the set of vertices in the CRT whose label is greater than $h$
(for this correspondence to hold, one needs to shift the labels so that the minimal label is $0$, and one
also re-roots the CRT at the vertex with minimal label). 
It follows that $\mathbf{N}_{h,\ve}$ counts those
connected components of the set of vertices with label greater than $h$ that contain 
(at least) one vertex with label $h+\ve$. In such a component, there is a unique 
vertex with label $h$ that is at minimal distance from the root, and, in the superprocess
setting, the remark following Theorem \ref{SBM} shows that this vertex corresponds to an upcrossing from $h$ to $h+\ve$. 

The paper is organized as follows. Section \ref{prelim} recalls basic facts about the Brownian snake,
which is our key tool to generate both the Brownian labels on the CRT and the historical paths of
super-Brownian motion. In Section \ref{sec-upcrossing}, we introduce upcrossings of the 
Brownian snake, and we state \cref{appro-LT}, which deals with the convergence \eqref{intro-2} under the excursion
measure of the Brownian snake. The proof of  \cref{appro-LT} is given in Section \ref{proof-appro},
after an important preliminary lemma (\cref{prelilemma}) has been established in Section \ref{prelimi}. Theorem \ref{SBM} is then an
easy consequence of \cref{appro-LT}.
Section \ref{sec-condi} provides conditional versions of \eqref{intro-2}, concerning first the 
excursion measure of the Brownian snake conditioned to have a fixed duration, and then
the same excursion measure under the additional conditioning that
the Brownian snake stays on the positive half-line. The latter conditional version
is needed for our application to the Brownian map in Section 
\ref{sec-applimap}, where we prove \cref{connec-compo}.

\section{Preliminaries about the Brownian snake}
\label{prelim}

We refer to the book \cite{LZurich} (especially Chapters IV and V) for the basic facts 
about the Brownian snake that we will use. 

\medskip
\noi{\bf The Brownian snake.} Throughout this work, $W=(W_s)_{s\geq 0}$ denotes the
one-dimensional Brownian snake. This is a strong Markov process taking values in the
space $\W$ of all finite continuous paths $\w:[0,\zeta]\la \R$, 
where $\zeta=\zeta_{(\w)}$ is a nonnegative real number depending on $\w$
and called the lifetime of $\w$. 
We write $\widehat \w:=\w(\zeta_{(\w)})$ for the
endpoint of $\w$. We let $(\zeta_s)_{s\geq 0}$ stand for the lifetime process
associated with $(W_s)_{s\geq 0}$, that is, $\zeta_s = \zeta_{(W_s)}$ for every $s\geq 0$.
For every $x\in\R$, we  identify the trivial element of $\W$ starting from $x$ and with zero lifetime with
the point $x$. 

It will be convenient to assume that the Brownian snake $(W_s)_{s\geq 0}$ is 
the canonical process on the space $C(\R_+,\W)$ of all continuous mappings
from $\R_+$ into $\W$. 
The notation $\P_x$ will then stand for the probability measure on $C(\R_+,\W)$ under which the
Brownian snake starts from $x$. Under $\P_x$, the process $(\zeta_s)_{s\geq 0}$ is a reflected Brownian motion on $\R_+$ started from $0$. Informally, the path $W_s$ is shortened from its tip when $\zeta_s$ decreases and, when
$\zeta_s$ increases, it
is extended by adding ``little pieces of Brownian paths'' at its tip. See \cite[Section IV.1]{LZurich} for a
more rigorous presentation.

We let $\N_x$ denote the (infinite) excursion measure of the
Brownian snake away from $x$. This excursion measure is normalized as in
\cite{LZurich}, so that, for every $\ve >0$,
$$\N_x\Big(\sup_{s\geq 0} \zeta_s > \ve\Big) = \frac{1}{2\ve}.$$
We also set 
$$\sigma=\inf\{s>0:\zeta_s=0\}$$
which represents the duration of the excursion under $\N_x$. 
The preceding informal description of the behavior of the Brownian snake remains valid under $\N_x$, but the ``law'' of the lifetime process under
$\N_x$ is now the It\^o measure of positive excursions of linear Brownian motion. 
Both under 
$\P_x$ and under $\N_x$, the Brownian snake takes values in the subset $\W_x$
of $\W$ that consists of all finite paths starting from $x$. Note that $W_s=x$ for every $s\geq \sigma$, $\N_x$-a.e.

For every $h\in\R$, we set
$$T_h:=\inf\{s\geq 0: \wh W_s =h\}$$
with the usual convention $\inf\varnothing=\infty$ that will be used throughout this work.  Suppose that $h\not = 0$. Then
\begin{equation}
\label{law-max}
\N_0(T_h<\infty) = \frac{3}{2 h^2}
\end{equation}
(see e.g. \cite[Lemma 2.1]{LGW}), and we will use  the notation $\N^h_0$
for the conditional probability measure
$$\N^h_0:=\N_0(\cdot\,|\,T_h<\infty).$$

\medskip
\noi{\bf Exit measures and the special Markov property.} We will make an extensive use of exit measures of the Brownian snake. Let $D$ be an open interval of $\R$, such that $D\not =\R$. Suppose that
$x\in D$ and, for every $\w\in \W_x$, set
$$\tau(\w)=\inf\{t\in [0,\zeta_{(\w)}]:\w(t)\notin D\},$$ 
where we recall that $\inf\varnothing=\infty$.
The exit measure $\z^D$ from $D$ (see \cite[Chapter 5]{LZurich}) is a random measure supported on $\partial D$, 
which is defined under $\N_x$ and is supported on the set of all exit points $W_s(\tau(W_s))$  for the paths $W_s$ such that  $\tau(W_s)<\infty$ (note that here $\partial D$ has at most two points, but the
preceding discussion remains valid for the $d$-dimensional Brownian snake
and an arbitrary subdomain $D$ of $\R^d$).

The first-moment formula for exit measures states that, for any nonnegative measurable 
function $g$ on $\partial D$,
\begin{equation}
\label{first-mo}
\N_x(\langle \z^D, g\rangle)= E_x[g(B_{\tau_D})]
\end{equation}
where, in the right-hand side, $B=(B_t)_{t\geq 0}$ is a linear Brownian motion starting from $x$ under the probability measure $P_x$, 
and $\tau_D:=\inf\{t\geq 0: B_t\notin D\}$. 

 We will
use the fact that, for every  $y\in\partial D$,
\begin{equation}
\label{exit-point}
\{\langle\z^D,\mathbf{1}_{\{y\}}\rangle >0\}= \{\exists s\geq 0 : \tau(W_s)<\infty\hbox{ and } W_s(\tau(W_s))=y\}\,,
\quad \N_x\hbox{-a.e.}
\end{equation}
It is immediate from the
support property of the exit measure that the set in the left-hand side is a subset of the set in the right-hand side. So, to get the equality in \eqref{exit-point}, it suffices to show that both sets
have the same $\N_x$-measure. However, using the connections 
between the Brownian snake and partial differential equations \cite[Chapters V,VI]{LZurich}, one verifies that the $\N_x$-measure of either set 
solves, as  a function of $x$, the differential equation $u''=4u^2$ in $D$ with boundary values $\infty$
at $y$, and $0$ at the other end of $D$ (at $\infty$ if $D$ is unbounded). Since this
boundary value problem has a unique nonnegative solution, the desired result follows.

A crucial ingredient of our study is the special Markov property of the Brownian snake
\cite{LG0}.
In order to state this property, we first observe that, $\N_x$-a.e., the set
$$\{s\geq 0: \tau(W_s)<\zeta_s\}$$
is open and thus can be written as a union of disjoint open intervals
$(a_i,b_i)$, $i\in I$, where $I$ may be empty. From the properties of the Brownian 
snake, it is easy to verify that, $\N_x$-a.e. for every $i\in I$ and every
$s\in(a_i,b_i)$, 
$$\tau(W_s)=\tau(W_{a_i})=\zeta_{a_i},$$
and more precisely all paths $W_s$, $s\in[a_i,b_i]$ coincide up to
their exit time from $D$. For every $i\in I$, we then define an element
$W^{(i)}$ of $C(\R_+,\W)$ by setting
$$W^{(i)}_s(t) := W_{(a_i+s)\wedge b_i}(\zeta_{a_i}+t),\quad \hbox{for }
0\leq t\leq \zeta_{(W^i_s)}:=\zeta_{(a_i+s)\wedge b_i}-\zeta_{a_i}.$$
Informally, the $W^{(i)}$'s represent the ``excursions'' of the Brownian snake outside $D$
(the word ``outside'' is a little misleading here, because although these excursions start from a point of $\partial D$, they will 
typically come back inside $D$).

We also need to introduce a $\sigma$-field that contains the information about the 
paths $W_s$ before they exit $D$. To this end, we set, for every $s\geq 0$,
$$\gamma^D_s:=\inf\{r\geq 0: \int_0^r \mathrm{d}u\,\mathbf{1}_{\{\tau(W_u)\geq \zeta_u\}} > s\},$$
and we let $\mathcal{E}^D$ be the $\sigma$-field generated by the
process $(W_{\gamma^D_s})_{s\geq 0}$ and the class of all sets 
that $\N_x$-negligible for every $x\in D$. The random measure $\z^D$ is measurable with respect to $\mathcal{E}^D$ (see 
\cite[Proposition 2.3]{LG0}).

We now state the special Markov property \cite[Theorem 2.4]{LG0}.

\begin{proposition}
\label{SMP}
Under $\N_x$, conditionally on $\mathcal{E}^D$, the point measure
$$\sum_{i\in I} \delta_{W^{(i)}}$$
is Poisson with intensity 
$$\int \z^D(\mathrm{d}y)\,\N_y.$$
\end{proposition}

Thanks to this proposition, we can consider each excursion $W^i$ again as a Brownian snake excursion
starting from a point of $\partial D$
and, if $D'$ is another domain containing $\partial D$, we can consider the ``subexcursions'' of $W^i$ 
outside $D'$, and so on. Repeated applications of this idea will play an important role in
what follows.

\medskip
\noi{\bf Local times.} We consider the total occupation measure $\mathcal{O}$
of the process $\wh W$, which is defined under $\N_x$
by the formula
$$\mathcal{O}(A): = \int_0^\sigma \mathrm{d}s\,\mathbf{1}_A(\wh W_s),$$
for any Borel subset $A$ of $\R$. 
The random measure $\mathcal{O}$  under $\N_0(\cdot\,|\,\sigma=1)$
is sometimes called one-dimensional ISE for integrated super-Brownian excursion (see
\cite{Al} and \cite[Section IV.6]{LZurich}). 

We will use the fact that $\mathcal{O}$ has $\N_x$-a.e. a continuous density $(L^a)_{a\in\R}$:
\begin{equation}
\label{occu}
\mathcal{O}(A)  =\int_{A} \mathrm{d}a\,L^a,
\end{equation}
for any Borel subset $A$ of $\R$. 
This can be derived from regularity properties of super-Brownian motion (see Section 1 and the references therein). 
Alternatively, we can use Theorem 2.1 in \cite{BMJ}, which gives the existence of a continuous density for $\mathcal{O}$
under $\N_0(\cdot\,|\,\sigma=1)$ (it is of course easy to get rid of the conditioning by $\sigma =1$ via a scaling argument). 

\section{Upcrossings of the Brownian snake}
\label{sec-upcrossing}

Consider the Brownian snake $(W_s)_{s\geq0}$ under $\N_x$ or under 
$\P_x$, for some fixed $x\in \R$. 

\begin{definition}
\label{defup}
Let $h\in\R$ and $\ve>0$. We say that $s\geq 0$ is an upcrossing time of
the Brownian snake from $h$ to $h+\ve$ if $\wh W_s=h$ and if there
exists $s'\in (s,\infty)$ such that $\wh W_{s'}=h+\ve$, $\zeta_r>\zeta_s$
for every $r\in(s,s']$, and $W_{s'}(t)>h$ for every $t\in(\zeta_s,\zeta_{s'}]$.
\end{definition}

The time $s'$ in the definition is in general not uniquely determined by $s$. However, there is a smallest
possible value of $s'$ such that the properties stated in the definition hold. In what follows,
we will always assume that $s'$ is chosen in this way, and we will say that $s'$
is associated with the upcrossing time $s$. 

\smallskip
\rem Obviously, a stopping time cannot be an upcrossing time of $W$.
On the other hand, it is easy to see that we can find a countable collection $(T_1,T_2,\ldots)$
of stopping times such that the set of all times $s'$ associated with upcrossing times 
from $h$ to $h+\ve$ is contained in $\{T_1,T_2,\ldots\}$. This remark will be useful
at the end of Section \ref{proof-appro}. 

\medskip
The reader may have  noticed that the preceding definition seems rather different
from the definition of an upcrossing time for a {\em function} $w:[0,T]\la \R$, which
was given in Section 1 (we might have considered upcrossing times of the function 
$s\la \wh W_s$, but this is not what we want!). To relate both definitions, we observe that, if $s$ is an upcrossing time of the Brownian snake
from $h$ to $h+\ve$,
and if $s'$ is the associated time, then $\zeta_s$ is an upcrossing time 
of the function $t\la W_{s'}(t)$ from $h$ to $h+\ve$. Definition \ref{defup} is more easily understood 
if we interpret the Brownian snake as a tree-indexed Brownian motion. Let us 
explain this in detail, as the relevant objects will also be useful later (see e.g.~\cite[Sections 3 and 4]{LGM} for a more detailed account
of the considerations that follow).

We
argue under $\N_x$, so that the lifetime process $(\zeta_s)_{s\geq 0}$
is just a single Brownian excursion. The tree  coded by 
$(\zeta_s)_{s\geq 0}$ is the quotient space $\t_\zeta:=[0,\sigma]\,/\!\sim$, where
the equivalence relation $\sim$ is defined by
$$s\sim s'\ \hbox{ if and only if } \zeta_s=\zeta_{s'}= \min_{r\in[s\wedge s',s\vee s']} \zeta_r.$$
We let $p_\zeta$ stand for
the canonical projection from $[0,\sigma]$ onto $\t_\zeta$, and equip $\t_\zeta$ with the metric $d_\zeta$
defined  by 
$$d_\zeta(p_\zeta(s),p_\zeta(s')):= \zeta_s + \zeta_{s'} - 2\,\min_{r\in[s\wedge s',s\vee s']} \zeta_r,$$
for every $s,s'\in[0,\sigma]$.
Then $\t_\zeta$ is a compact $\R$-tree, which is rooted at $\rho_\zeta:=p_\zeta(0)$. 
Note that the generation (distance from the root) of the vertex $p_\zeta(s)$ is $\zeta_s$. For $a,b\in\t_\zeta$,
we will use the notation $\llbracket a,b \rrbracket$ for the line segment between $a$ and $b$ in $\t_\zeta$.
The notions of an 
ancestor and a descendant in $\t_\zeta$ are defined in an obvious way: For $a,b\in\t_\zeta$, $a$ is
an ancestor of $b$ if $a$ belongs to $\llbracket \rho_\zeta,b \rrbracket$. 
If $s,s'\in[0,\sigma]$, $p_\zeta(s)$ is an ancestor of $p_\zeta(s')$ if and only if
$\zeta_r\geq \zeta_s$ for every $r\in[s\wedge s',s\vee s']$.

It follows from the properties of the Brownian snake that, $\N_x$-a.e., 
$\wh W_s=\wh W_{s'}$ for every $s,s'$ such that $s\sim s'$.
Hence we can define $\Gamma_a$ for every $a\in\t_\zeta$ by declaring
that $\Gamma_{p_\zeta(s)}= \wh W_s$ for every $s\in[0,\sigma]$, and it is very
natural to interpret $(\Gamma_a)_{a\in\t_\zeta}$ as Brownian motion
indexed by $\t_\zeta$. We vew $\Gamma_a$ as a spatial location or label assigned 
to the vertex $a$. For every $s\in[0,\sigma]$ and every $t\in[0,\zeta_s]$, $W_s(t)$ corresponds to
the spatial location of the ancestor of $p_\zeta(s)$ at generation $t$.

It is now easy to verify that upcrossing times of $W$ from $h$ to $h+\ve$
are in one-to-one correspondence with vertices $a$ of $\t_\zeta$
such that $\Gamma_a=h$ and there exists a descendant $b$ of 
$a$ in $\t_\zeta$ such that $\Gamma_{b}=h+\ve$ and $\Gamma_c>h$
for every interior point $c$ of the line segment $\llbracket a,b \rrbracket$. 
In this form, we see that our definition is the exact analog
of the one for upcrossing times of a real function defined on the interval $[0,T]$ (provided
we see $[0,T]$ as an $\R$-tree rooted at $0$). 

\begin{lemma}
\label{fini-upcro}
Let $h\in\R$ and $\ve>0$. Let $N_{h,\ve}$ be the number of upcrossing times
of the Brownian snake from $h$ to $h+\ve$. Then, $N_{h,\ve}<\infty$, $\N_x$-a.e.
\end{lemma}

\proof By continuity, there exists $\N_x$-a.e.~a real $\delta>0$ such that
$|\wh W_{s_1} - \wh W_{s_2}|<\ve$ for every $s_1,s_2\geq 0$ such that $|s_1-s_2|\leq\delta$.
If $s$ is an upcrossing time from $h$ to $h+\ve$, we let $s'>s$ be the time associated with $s$ (see the comment
following Definition \ref{defup}), and we set $I_s:=[s'-\delta,s'+\delta]$. The statement of the lemma follows
from the fact that the intervals $I_s$,
when $s$ varies over the set of all upcrossing times, are pairwise disjoint.
To verify the latter fact, consider two upcrossing times $s$ and $\tilde s$
and the associated times $s'$ and $\tilde s'$. If $s<s'<\tilde s<\tilde s'$, the desired
property is immediate from our choice of $\delta$ and the definition of upcrossing times.
From this definition, it is also easy to
verify that we cannot have $s<\tilde s<s'<\tilde s'$ (otherwise, $p_\zeta(\tilde s)$ would be
both a strict descendant of $p_\zeta(s)$ and an ancestor of $p_\zeta(s')$, implying $\wh W_{\tilde s}>h$).
The only case that remains is when $s<\tilde s<\tilde s'< s'$. In that case,
$p_\zeta(\tilde s)$ is an ancestor of $p_\zeta(\tilde s')$ but not an ancestor of
$p_\zeta(s')$. It follows that, if $s'':=\inf\{r\geq \tilde s: \zeta_r<\zeta_{\tilde s}\}$, we have
$\tilde s<\tilde s'<s''<s'$. However, $p_\zeta(s'')=p_\zeta(\tilde s)$ and so $\wh W_{s''}=\wh W_{\tilde s}=h$,
whereas $\wh W_{\tilde s'}= h+\ve$ and $\wh W_{s'}=h+\ve$. The property
$s'-\tilde s'>2\delta$ now follows from our choice of $\delta$.
\endproof

Recall our notation $\N_0^h$ for the excursion measure $\N_0$ conditioned on the
event that the Brownian snake hits the level $h$. The following statement is the main 
technical result of the paper, from which we will deduce the theorems stated in Section 1.

\begin{theorem}
\label{appro-LT}
Let $h\in\R\backslash\{0\}$. We have
$$\ve^3\,N_{h,\ve} \;\build{\la}_{\ve \to 0}^{} \;\frac{c_1}{2}\,L^h,$$
in probability under $\N^h_0$. 
Here $L^h$ is the density at $h$ of the occupation measure $\mathcal{O}$, and the constant $c_1>0$ is determined by the identity
$$\int_0^\infty \frac{\mathrm{d}u}{\sqrt{\frac{8}{3}u^3 + (c_1)^2}} = 1.$$
\end{theorem}

\rem We exclude the value $h=0$, in particular because the measure $\N^h_0$ is not defined when $h=0$. 

\smallskip

The proof of  
\cref{appro-LT} is given below in Section \ref{proof-appro}. Section \ref{prelimi}
contains some preliminary lemmas.

\section{Preliminary lemmas}
\label{prelimi}

For technical reasons, we will first deal with the Brownian
snake under the probability measure $\P_0$. We write $(\ell^0_s)_{s\geq 0}$ for the local time at 
level $0$ of the reflected Brownian motion $(\zeta_s)_{s\geq 0}$ (the normalization of local times is such that
the occupation density formula holds, and local times are right-continuous in the space variable). For every
$r>0$, we set
$$\eta_r:=\inf\{s\geq 0: \ell^0_s > r\}.$$
The excursions of $W$ away from $0$, before time $\eta_r$, form a Poisson measure with 
intensity $r\,\N_0$. 

For every $\ve>0$ and $\w\in \W$, set $\tau_\ve(\w):=\inf\{t\in[0,\zeta_{(\w)}]: \w(t)\geq \ve\}$, and also define, for every
$r>0$,
$$M_\ve(r):=\#\{s\in[0,\eta_r]: \tau_\ve(W_s)=\infty\hbox{ and }s\hbox{ is an upcrossing time of }W\hbox{ from }0\hbox{ to }\ve\}.$$
Lemma \ref{fini-upcro} implies that $M_{\ve}(r)<\infty$, $\P_0$-a.s. (note that only finitely many 
excursions of $W$ away from $0$ hit $\ve$ before time $\eta_r$).

\begin{lemma}
\label{prelilemma}
For every $\ve>0$ and $r>0$,
$$\E_0[ M_\ve(r)] = c_1\,r\,\ve^{-2}\,,$$
where the constant $c_1$ is as in \cref{appro-LT}.
\end{lemma}

\proof In this proof, $\ve>0$ and $r>0$ are fixed. We also consider
a real $\delta>0$, that later will tend to $0$ (to avoid problems with uncountable
unions of negligible sets, we may and will restrict our attention to
rational values of $\delta$). We write
$$\n^0(d\omega):= \sum_{i\in I_0} \delta_{\omega^0_i}(d\omega)$$
for the point measure of excursions of $W$ away from $0$ before time
$\eta_r$. With each excursion $\omega^0_i$, we associate its exit
measure $\z^{(-\delta,\ve)}(\omega^0_i)$ from the interval
$(-\delta,\ve)$. This exit measure is a finite measure supported on the pair $\{-\delta,\ve\}$. We set
$$X^1_\delta :=\sum_{i\in I_0} \langle\z^{(-\delta,\ve)}(\omega^0_i), {\mathbf 1}_{\{-\delta\}}\rangle,$$
which represents the total mass assigned to the point $-\delta$ by the exit measures
associated with the excursions $\omega^0_i$, $i\in I_0$. 

Then, for every $i\in I_0$ (we need only consider those values of $i$ such that $\omega^0_i$ hits $-\delta$), we can
introduce the excursions of $\omega^0_i$ outside $(-\delta,\ve)$ that start from $-\delta$, as defined in Section 
\ref{prelim}. Write $(\tilde\omega^0_j)_{j\in J_0}$ for the collection of all these excursions 
when $i$ varies over $I_0$. By the special Markov property (\cref{SMP}), we know that, conditionally
on $X^1_\delta$, the point measure
$$\sum_{j\in J_0} \delta_{\tilde\omega^0_j}$$
is Poisson with intensity $X^1_\delta\,\N_{-\delta}$. 

For every $j\in J_0$, $\wt\omega^0_j$ is a Brownian snake excursion starting from $-\delta$, and 
therefore we can consider its exit measure $\z^{(-\infty,0)}(\tilde\omega^0_j)$ from the
interval $(-\infty,0)$. We then set
$$Y^1_\delta:= \sum_{j\in J_0} \langle\z^{(-\infty,0)}(\tilde\omega^0_j) , 1\rangle.$$
Furthermore, for every $j\in J_0$, we can also consider the excursions of $\tilde\omega^0_j$
outside $(-\infty,0)$ (of course these excursions start from $0$). We write $(\omega^1_i)_{i\in I_1}$
for the collection of all these excursions when $j$ varies over $J_0$, and we set
$$\n^1_\delta(d\omega)=\sum_{i\in I_1} \delta_{\omega^1_i}(d\omega).$$
By the special Markov property again, we get that, conditionally on $Y^1_\delta$, the point
measure $\n^1_\delta$ is Poisson with intensity $Y^1_\delta\,\N_0$. 
Informally, the point measure $\n^1_\delta$ contains the information about the behavior after
their first return to $0$ via $-\delta$ of those paths $W_s$ that hit $-\delta$ before hitting $\ve$.

We can continue this construction by induction. Let us briefly describe the second step. We set
$$X^2_\delta := \sum_{i\in I_1} \langle\z^{(-\delta,\ve)}(\omega^1_i), {\mathbf 1}_{\{-\delta\}}\rangle,$$
and write $(\tilde\omega^1_j)_{j\in J_1}$ for the collection of all excursions of $\omega^1_i$, $i\in I_1$,
outside $(-\delta,\ve)$ that start from $-\delta$. We then set
$$Y^2_\delta:= \sum_{j\in J_1} \langle\z^{(-\infty,0)}(\tilde\omega^0_j) ,1\rangle,$$
and 
$$\n^2_\delta(d\omega)=\sum_{i\in I_2} \delta_{\omega^2_i}(d\omega)$$
where $(\omega^2_i)_{i\in I_2}$ is the collection of all excursions of $\tilde\omega_j$, $j\in J_1$,
outside $(-\infty,0)$. Again, conditionally on $Y^2_\delta$, the point measure
$\n^2_\delta$ is Poisson with intensity $Y^2_\delta\,\N_0$. 

At every step $k\geq 1$, we similarly get a nonnegative random variable $Y^k_\delta$, and 
a point measure 
$$\n^k_\delta(d\omega)=\sum_{i\in I_k} \delta_{\omega^k_i}(d\omega),$$
which, conditionally on $Y^k_\delta$, is Poisson with intensity $Y^k_\delta\,\N_0$. Informally,
$\n^k_\delta$ describes the paths $W_s$ after their $k$-th return to $0$ via $-\delta$, for those
paths $W_s$ that perform $k$ descents from $0$ to $-\delta$ before they hit $\ve$. 

We now set, for every integer $k\geq 0$,
$$M^k_{\ve,\delta} := \sum_{i\in I_k} {\mathbf 1}_{\{\langle\z^{(-\delta,\ve)}(\omega^k_i), {\mathbf 1}_{\{\ve\}}\rangle
>0\}},$$
which counts those Brownian snake excursions $\omega^k_i$, $i\in I_k$, for which there exists $s\geq 0$
such that the path $\omega^k_i(s)$ hits $\ve$ before $-\delta$ (by \eqref{exit-point}, the existence of such a value of $s$
is equivalent to the property $\langle\z^{(-\delta,\ve)}(\omega^k_i), {\mathbf 1}_{\{\ve\}}\rangle>0$).
We also set
$$M_{\ve,\delta}:=\sum_{k=0}^\infty M^k_{\ve,\delta}.$$
At this point, we need another lemma.

\begin{lemma}
\label{prelilemma2}
We have $M_{\ve,\delta}\leq M_\ve(r)$ for every $\delta>0$. Moreover,
$$M_\ve(r)=\lim_{\delta\downarrow 0} M_{\ve,\delta},\quad\P_0\hbox{ a.s.}$$
\end{lemma}

We postpone the proof of \cref{prelilemma2} and complete the proof
of \cref{prelilemma}. We note that, by Lemma \ref{prelilemma2} and Fatou's lemma,
we have
$$\E_0[M_\ve(r)]\leq \liminf_{\delta\downarrow 0} \E_0[M_{\ve,\delta}].$$
On the other hand, the first assertion of  \cref{prelilemma2} also shows that
$\E_0[M_{\ve,\delta}]\leq \E_0[M_\ve(r)]$ for every $\delta>0$, so that we have
\begin{equation}
\label{prelitech1}
\E_0[M_\ve(r)]= \lim_{\delta\downarrow 0} \E_0[M_{\ve,\delta}].
\end{equation}
\par To complete the argument, we will compute $\E_0[M_{\ve,\delta}]$. 
We first set
$$a_{\ve,\delta}:=\N_0(\langle\z^{(-\delta,\ve)},\mathbf{1}_{\{\ve\}}\rangle >0).$$
As we already mentioned after \eqref{exit-point}, we have $a_{\ve,\delta}=u(0)$, where the function
$(u(x),x\in(-\delta,\ve))$ solves the differential equation $u''=4\,u^2$ with boundary
conditions $u(\ve)=\infty$ and $u(-\delta)=0$. Solving this differential equation 
leads to
$$\int_0^{u(x)} \frac{\mathrm{d}u}{\sqrt{\frac{8}{3}u^3 +  c_{\ve,\delta}}} = x + \delta\;,\quad x\in(-\delta,\ve),$$
where the constant $c_{\ve,\delta}>0$ is determined by
$$\int_0^{\infty} \frac{\mathrm{d}u}{\sqrt{\frac{8}{3}u^3 +  c_{\ve,\delta}}} = \ve + \delta.$$
It follows that $c_{\ve,\delta}= (\ve +\delta)^{-6} (c_1)^2$, where $c_1$ is as in the
statement of \cref{appro-LT}. Since $a_{\ve,\delta}=u(0)$, we have then
$$\int_0^{a_{\ve,\delta}} \frac{\mathrm{d}u}{\sqrt{\frac{8}{3}u^3 + (\ve +\delta)^{-6} (c_1)^2}} = \delta,$$
and elementary analysis shows that
\begin{equation}
\label{prelitech2}
\lim_{\delta\to0} \delta^{-1}\,a_{\ve,\delta} = \ve^{-3}\,c_1.
\end{equation}

Now note that
$$\E_0[M^0_{\ve,\delta}] = a_{\ve,\delta}\,r$$
and, using the conditional distribution of $\n^k_\delta$ given $Y^k_\delta$,
$$\E_0[M^k_{\ve,\delta}]= a_{\ve,\delta}\,\E_0[Y^k_\delta],$$
for every $k\geq 1$. On the other hand, by the first moment formula for exit 
measures \eqref{first-mo}, we have
$$\E_0[X^1_\delta]= \frac{\ve}{\ve+\delta}\,r$$
and 
$$\E_0[Y^1_\delta]= \E_0[X^1_\delta]= \frac{\ve}{\ve+\delta}\,r.$$
An easy induction argument gives, for every $k\geq1$,
$$\E_0[Y^k_\delta]= \Big(\frac{\ve}{\ve+\delta}\Big)^k\,r.$$
Hence,
$$\E_0[M_{\ve,\delta}]=
\sum_{k=0}^\infty \E_0[M^k_{\ve,\delta}]= a_{\ve,\delta} \sum_{k=0}^\infty \Big(\frac{\ve}{\ve+\delta}\Big)^k\,r
= a_{\ve,\delta}\,r\Big(\frac{\ve+\delta}{\delta}\Big).$$
Using \eqref{prelitech1} and \eqref{prelitech2}, we get
$$\E_0[M_\ve(r)]=\lim_{\delta\downarrow0} \E_0[M_{\ve,\delta}] = c_1r\,\ve^{-2},$$
which completes the proof of Lemma \ref{prelilemma}. \endproof

\medskip
\noindent{\it Proof of Lemma \ref{prelilemma2}.} Recall the construction of excursions 
of the Brownian snake outside
an interval. For every $k\geq 0$ and every $i\in I_k$, the excursion $\omega^k_i$
corresponds to a closed subinterval $\mathcal{I}_{k,i}$ of $[0,\eta_r]$ (in such a way
that the paths $\omega^k_i(s)$, $s\geq 0$, are exactly the paths $W_s$, $s\in\mathcal{I}_{k,i}$
shifted at the time of their $k$-th return to $0$ via $-\delta$). Next, if 
$\langle\z^{(-\delta,\ve)}(\omega^k_i), {\mathbf 1}_{\{\ve\}}\rangle>0$, we can find $s_0\in\mathcal{I}_{k,i}$
such that $\wh W_{s_0}=\ve$ and $\tau_\ve(W_{s_0})= \zeta_{s_0}$, and the path 
$W_{s_0}$ performs exactly $k$ descents from $0$ to $-\delta$. Set
$$\lambda_0(W_{s_0})=\sup\{t\in[0,\zeta_{s_0}]: W_{s_0}(t)=0\},$$
and
$$r_0=\sup\{s\leq s_0:\zeta_s=\lambda_0(W_{s_0})\}.$$
Note that $r_0$ also belongs to $\mathcal{I}_{k,i}$, because, for $r_0 < s\leq s_0$, we have
$\zeta_s>\zeta_{r_0}=\lambda_0(W_{s_0})$ and the path
$W_s$ coincides with $W_{s_0}$ up to a time strictly greater than $\lambda_0(W_{s_0})$. From our definitions,
$r_0$ is an upcrossing time of $W$ from $0$ to $\ve$. If we now vary $k$ and $i$
(among all pairs $(k,i)$ such that $\langle\z^{(-\delta,\ve)}(\omega^k_i), {\mathbf 1}_{\{\ve\}}\rangle>0$),
we get distinct upcrossing times. This is obvious if we vary $i$ for a fixed value of $k$, because 
the intervals $\mathcal{I}_{k,i}$, $i\in I_k$, are disjoint. If we vary $k$, this follows from the fact that
$k$ can be interpreted as the number of descents of $W_{r_0}$ from $0$ to $-\delta$. The preceding
discussion shows that $M_\ve(r)\geq M_{\ve,\delta}$, proving the first assertion of the lemma.

In order to prove the second assertion, let us start with a few remarks. Suppose that
$s$ is an upcrossing time from $0$ to $\ve$, and let $s'$ be associated with
$s$ as explained after \cref{defup}. Let $k$ be the number of
descents from $0$ to $-\delta$ of the path $W_s$. Then $s$ must
belong to exactly one interval $\mathcal{I}_{k,i}$, with $i\in I_k$, and $s'$
belongs to the same interval. Write
$\mathcal{I}_{k,i}=[\alpha_{k,i},\beta_{k,i}]$, with $\alpha_{k,i}<\beta_{k,i}$,
and note that, by construction, all paths $W_u$, $u\in\mathcal{I}_{k,i}$ coincide
up to time $\zeta_{\alpha_{k,i}}=\zeta_{\beta_{k,i}}$. Furthermore, for 
every $u\in [0,\beta_{k,i}-\alpha_{k,i}]$, the path $\omega^k_i(u)$ is just the
path  $W_{\alpha_{k,i}+u}$ shifted at time $\zeta_{\alpha_{k,i}}$. Now, using the fact
that $W_{s'}$ makes the same number of descents from $0$ to $-\delta$
as $W_s$, and the definition of an upcrossing time, we see that $\omega^k_i(s'-\alpha_{k,i})$ hits $\ve$
before $-\delta$. Using \eqref{exit-point}, it follows that
\begin{equation}
\label{prelitech3}
\langle\z^{(-\delta,\ve)}(\omega^k_i), {\mathbf 1}_{\{\ve\}}\rangle >0.
\end{equation}

Consider then another upcrossing time $\tilde s> s$ and the associated time
$\tilde s'$.  To simplify notation, set
$$\check\zeta_{s,\tilde s} :=\min_{s\leq u\leq \tilde s} \zeta_u.$$
By the properties of the Brownian snake, the paths $W_s$ and $W_{\tilde s}$
coincide over the interval $[0, \check\zeta_{s,\tilde s}]$.  Suppose that $W_{\tilde s}$
also makes $k$ descents from $0$ to $-\delta$, and belongs to the same interval $\mathcal{I}_{k,i}$
as $s$. Then necessarily $\check\zeta_{s,\tilde s}\geq \zeta_{\alpha_{i,k}}$, and we have
\begin{equation}
\label{prelitech4}
\min\{W_{\tilde s}(t): t\in [\check\zeta_{s,\tilde s},\zeta_{\tilde s}]\} > -\delta,
\end{equation}
because otherwise $W_{\tilde s}$ would make (at least) $k+1$ descents from $0$ to $-\delta$. 

Now let $s_1,s_2,\ldots,s_p$ be $p$ distinct upcrossing
times from $0$ to $\ve$ such that $s_1<s_2<\cdots<s_p$. The second assertion of the lemma will follow
if we can prove that we have $M_{\ve,\delta} \geq p$ for $\delta>0$ small enough.  We first observe that, for
every $i,j\in\{1,\ldots,p\}$ such that $i<j$, we have $\check\zeta_{s_i,s_j}<\zeta_{s_j}$, because
otherwise (recalling the definition of an upcrossing time) $s_j$ would be a time
of local minimum of $\zeta$, and it is easy to see that such a time cannot be an upcrossing time. Then it also
follows that we have
$$\min\{W_{s_j}(t): t\in [\check\zeta_{s_i,s_j},\zeta_{s_j}]\}<0.$$
Indeed argue by contradiction and suppose that the latter minimum vanishes. Then 
writing $s'_j$ for the time associated with $s_j$, we obtain that the path $W_{s'_j}$ has a local
minimum  equal to $0$ at time $\zeta_{s_j}<\zeta_{s_j'}$. This is a contradiction because, with 
probability one, none of the paths $W_s$ can have a local minimum equal to $0$ at an interior point of $[0,\zeta_s]$
(observe that it is enough to consider rational values of $s$, and then note that a fixed constant is 
a.s.~not a local minimum of linear Brownian motion).

 To complete the argument, we observe that, if 
 $$-\delta > \max_{1\leq i<j\leq p} \Big(\min\{W_{s_j}(t): t\in [\check\zeta_{s_i,s_j},\zeta_{s_j}]\}\Big),$$
 then the pairs $(k_j,i_j)$ corresponding to the different upcrossing times $s_1,\ldots,s_p$
 must be distinct, because otherwise this would contradict the property \eqref{prelitech4}. Furthermore,
 we can apply \eqref{prelitech3} to each pair $(k_j,i_j)$, and it follows that, for $\delta>0$ small enough, we 
 have $M_{\ve,\delta} \geq p$. This completes the proof of Lemma \ref{prelilemma2}. \hfill$\square$

\section{Proofs of Theorem
\ref{appro-LT} and Theorem \ref{SBM}}
\label{proof-appro}

Most of this section is devoted to the proof of \cref{appro-LT}. We then 
explain how to derive Theorem \ref{SBM} from this statement.

\medskip

\noi{\it Proof of 
\cref{appro-LT}.}
Let $M_\ve$ be the analog of $M_\ve(r)$ for a single Brownian snake excursion,
$$M_\ve:= \#\{s\in[0,\sigma]: \tau_\ve(W_s)=\infty\hbox{ and }s\hbox{ is an upcrossing time of }W\hbox{ from }0\hbox{ to }\ve\}.$$
Clearly, $\E_0[M_\ve(r)]= r\,\N_0(M_\ve)$, and we deduce from Lemma \ref{prelilemma} that
\begin{equation}
\label{appro-tech0}
\N_0(M_\ve)= c_1\,\ve^{-2}.
\end{equation}

Let us fix $h\in \R\backslash\{0\}$. In the present  proof, we argue under the 
probability measure $\N_0^{h}$ (for technical reasons it will be convenient to enlarge the 
probability space so that it carries real random variables independent of the Brownian snake).  We have 
$$N_{h,\ve}= N^0_{h,\ve} + N^1_{h,\ve}+ \cdots + N^k_{h,\ve} + \cdots$$
where, for every integer $k\geq 0$, $N^k_{h,\ve}$ counts the number of upcrossing
times $s$ from $h$ to $h+\ve$ such that the path $W_s$ has made exactly 
$k$ upcrossings from $h$ to $h+\ve$. Considering the excursions of the Brownian snake outside the domain $(-\infty,h)$ if 
$h>0$, or the domain $(h,\infty)$ if $h<0$, we see that
$$N^0_{h,\ve}=N^{0,1}_{h,\ve} + N^{0,2}_{h,\ve} + \cdots + N^{0,n^0_\ve}_{h,\ve},$$
where $n^0_\ve$ denotes the number of excursions outside $(-\infty,h)$ (if $h>0$)
or outside $(h,\infty)$ (if $h<0$) that hit $h+\ve$, and, for every $1\leq i\leq n^0_\ve$,
$N^{0,i}_{h,\ve}$ counts the contribution to $N^0_{h,\ve}$ of the $i$-th excursion, assuming that these excursions are
listed in a uniform random order\footnote{The argument also goes through if excursions 
are listed in chronological order, but then we need a slightly more precise version of the
special Markov property.} given the Brownian snake.  In other words,
$N^{0,i}_{h,\ve}$ counts those upcrossing times $s$ that belong to the interval associated
with the $i$-th excursion, and have the additional property that $W_s$ makes no upcrossing
from $h$ to $h+\ve$. From the special Markov property, we
see that, conditionally on $n^0_\ve$, the variables $N^{0,1}_{h,\ve}, N^{0,2}_{h,\ve}, \ldots$
are independent and follow the distribution of $M_\ve$ under $\N_0(\cdot\mid T_\ve<\infty)$.
By scaling, the latter distribution does not depend on $\ve$, and we denote it by $\mu$.

A similar decomposition holds for $N^k_{h,\ve}$, for every 
$k\geq 1$. Let us discuss the case $k=1$. We consider again
the excursions of the Brownian snake outside the domain $(-\infty,h)$ if 
$h>0$, or the domain $(h,\infty)$ if $h<0$. We apply the special Markov property to
each of these excursions (which start from $h$) and to the domain 
$(-\infty,h+\ve)$, in order to get a collection of Brownian snake excursions 
starting from $h+\ve$. Once again, we apply the special Markov property
to each of the latter excursions and to the domain $(h,\infty)$, and we let
$n^1_\ve$ be the number of the resulting excursions (starting from $h$) that hit $h+\ve$.
We have then
$$N^1_{h,\ve}=N^{1,1}_{h,\ve} + N^{1,2}_{h,\ve} + \cdots + N^{1,n^1_\ve}_{h,\ve},$$
where, conditionally on the pair $(n^0_\ve,n^1_\ve)$, the random variables
$N^{1,1}_{h,\ve},N^{1,2}_{h,\ve},\ldots$ are independent and distributed according to $\mu$.  
Furthermore, still conditionally on $(n^0_\ve,n^1_\ve)$, the vector
$(N^{1,1}_{h,\ve},\ldots ,N^{1,n^1_\ve}_{h,\ve})$ is independent of the
vector $(N^{0,1}_{h,\ve},\ldots ,N^{0,n^0_\ve}_{h,\ve})$: This follows again from the
special Markov property, and the fact that the
vector $(N^{0,1}_{h,\ve},\ldots ,N^{0,n^0_\ve}_{h,\ve})$ is measurable with respect to 
the $\sigma$-field generated by the paths $W_s$ up to the end of their first 
upcrossing from $h$ to $h+\ve$. 

Arguing inductively, we get that, for every $k\geq 0$,
\begin{equation}
\label{appro-tech1}
N^k_{h,\ve} = \sum_{i=1}^{n^k_\ve} N^{k,i}_{h,\ve}
\end{equation}
where, conditionally on the sequence $(n^0_\ve,n^1_\ve,\ldots,n^k_\ve)$, the random variables
$N^{k,i}_{h,\ve}$,  $1\leq i\leq n^k_\ve$, are independent, and independent of the 
collection $(N^{\ell,j}_{h,\ve})_{0\leq \ell<k,1\leq j\leq n^\ell_\ve}$, and are distributed 
according to $\mu$.
Furthermore, the variables $n^k_\ve$ can be characterized as follows. For every
$k\geq 0$, $n^k_\ve$ counts the instants $s$ such that:
\begin{enumerate}
\item[(i)] the path $W_s$ makes exactly $k$ upcrossings from $h$ to $h+\ve$;
\item[(ii)] $\zeta_s$ is the time of the first return of $W_s$ to $h$ after the $k$-th upcrossing
from $h$ to $h+\ve$
(when $k=0$, $\zeta_s$ coincides the first hitting time of $h$ by $W_s$);
\item[(iii)] in the tree $\t_\zeta$, $p_\zeta(s)$ has (at least) one descendant $p_\zeta(s')$ such that 
$\wh W_{s'}= h+\ve$.
\end{enumerate}
Notice that, for every fixed $\ve>0$, we have $n^k_\ve=0$ for $k$ large enough, $\N^h_0$-a.s.

Let us emphasize that there is no independence 
between the sequence $(n^k_\ve)_{k\geq 0}$ on one hand and the collection of variables $(N^{k,i}_{h,\ve})_{k\geq0,1\leq i\leq n^k_\ve}$
on the other hand. At an intuitive level, if $N^{0,1}_{h,\ve}$ is large,  the exit measure from $(-\infty,h+\ve)$ of the first excursion
outside $(-\infty,h)$ (if $h>0$) or outside $(h,\infty)$ (if $h<0$)
is likely to be large, and, with high probability, $n^1_\ve$ will also be large. 

\begin{lemma}
\label{key-approx}
For every $\ve>0$, set
$$n_\ve:=\sum_{k=0}^\infty n^k_\ve\,.$$
We have
$$\frac{4\ve^3}{3}\,n_\ve \build{\la}_{\ve \to0}^{} L^h$$
in probability under $\N^{h}_0$.
\end{lemma}

We postpone the proof of the lemma and complete the proof of \cref{appro-LT}.
Write $\xi^\ve_1,\xi^\ve_2,\ldots$ for the sequence 
$$N^{0,1}_{h,\ve},N^{0,2}_{h,\ve},\ldots, N^{0,n^0_\ve}_{h,\ve},
N^{1,1}_{h,\ve},N^{1,2}_{h,\ve},\ldots, N^{1,n^1_\ve}_{h,\ve},N^{2,1}_{h,\ve},N^{2,2}_{h,\ve},\ldots, N^{2,n^2_\ve}_{h,\ve},\ldots$$
which is completed by adding a sequence of independent random variables distributed according to $\mu$
at its end (these auxiliary random variables are supposed to be independent of the Brownian snake). 
As a consequence of the properties stated after \eqref{appro-tech1}, it is a simple exercise to
verify that $\xi^\ve_1,\xi^\ve_2,\ldots$ form a sequence of independent random variables distributed
according to $\mu$. Set $S^\ve_j=\xi^\ve_1+\xi^\ve_2+\cdots+ \xi^\ve_j$, for every $j\geq 0$. By construction,
we have 
$$N_{h,\ve} = S^\ve_{n_\ve}.$$
To complete the proof, we simply use the law of large numbers. Using \eqref{law-max} and 
\eqref{appro-tech0}, we get that the first moment of $\mu$ is
$$\N_0(M_\ve\mid T_\ve<\infty)=\frac{\N_0(M_\ve)}{\N_0(T_\ve<\infty) }= \frac{2c_1}{3}$$
which does not depend on $\ve$ as expected. It is easy to verify that $n_\ve \la \infty$
 as $\ve \to0$, $\N^{h}_0$ a.s.~(by \eqref{exit-point} and the special Markov property, this is even true
 if we replace $n_\ve$ by $n^0_\ve$), and the law of large numbers gives
 $$\frac{1}{n_\ve}\,N_{h,\ve} = \frac{1}{n_\ve}\,S^\ve_{n_\ve}\build{\la}_{\ve\to0}^{} \frac{2c_1}{3}$$
 in $\N^{h}_0$-probability. Notice that the preceding argument applies even though
 $n_\ve$ is {\em not} independent of the sequence $(S^\ve_j)_{j\geq 1}$. The convergence of 
 \cref{appro-LT} now follows by writing
 $$\ve^{3} N_{h,\ve} = (\ve^3 n_\ve)\times \frac{1}{n_\ve}\,N_{h,\ve}$$
 and using Lemma \ref{key-approx}. \hfill$\square$
 
 \medskip
 \noindent{\it Proof of Lemma \ref{key-approx}.} In this proof, we assume that $h>0$. Only minor modifications
 are needed when $h<0$. It will be convenient to replace the convergence of Lemma \ref{key-approx} by
 an analogous convergence in terms of certain exit measures. We argue in a way very similar to the proof
 of Lemma \ref{prelilemma}, and for this reason we will omit some details. We first set
\begin{align*}
Z^{\ve,0}&:= \langle \z^{(-\infty,h)}, 1\rangle,\\
 Z^{\ve,1}&:= \langle \z^{(-\infty,h+\ve)}, 1\rangle.
 \end{align*}
 Then we consider all excursions of the Brownian snake outside $(-\infty,h+\ve)$ and we define
 $Z^{\ve,2}$ as the sum, over all these excursions, of the total masses of their exit measures
 from $(h,\infty)$. For each of the preceding excursions, we consider its ``subexcursions'' 
 outside $(h,\infty)$ and define $Z^{\ve,3}$ as the sum over all these (and over all the initial excursions outside 
 $(-\infty,h+\ve)$) of the total masses of their exit measures from $(-\infty,h+\ve)$. We continue by induction in an obvious way.
 Informally, for any $k\geq 1$,  $Z^{\ve,2k-1}$ ``counts'' the paths $W_s$ that make exactly 
 $k$ upcrossings from $h$ to $h+\ve$, and are stopped at the end of the $k$-th upcrossing, and
 similarly $Z^{\ve,2k}$ ``counts'' the paths $W_s$ that make $k$ upcrossings from $h$ to $h+\ve$,
 then one additional descent from $h+\ve$ to $h$, and are stopped at the end of this last descent. 
 
 Using a symmetry argument analogous to the classical reflection principle for Brownian motion
 (but now relying on the special Markov property rather than on the strong Markov property of Brownian motion),
 one immediately verifies that
 $$(Z^{\ve,0},Z^{\ve,1},Z^{\ve,2}) \build{=}_{}^{\rm(d)} ( \langle \z^{(-\infty,h)}, 1\rangle, 
  \langle \z^{(-\infty,h+\ve)}, 1\rangle,  \langle \z^{(-\infty,h+2\ve)}, 1\rangle).$$
  This argument is easily extended to yield
  $$(Z^{\ve,k})_{k\geq 0} \build{=}_{}^{\rm(d)} ( \langle \z^{(-\infty,h+k\ve)}, 1\rangle)_{k\geq 0}.$$
As an easy consequence of the special Markov property and the first-moment formula \eqref{first-mo}, the process $(\langle \z^{(-\infty,h+a)}, 1\rangle)_{a\geq0}$, which is
  now indexed by the real variable $a\geq 0$, is a nonnegative martingale (in fact a critical continuous-state branching process) under
  $\N^{h}_0$. Consequently, this process has a c\`adl\`ag modification, which we consider from now on, and using the preceding identity in distribution, we have, for every $\ve>0$,
  \begin{equation}
  \label{key-tech1}
  \sup_{k\geq 0} Z^{\ve,k}\; \build{\leq}_{}^{\rm(d)} \;\sup_{a\geq 0} \;\langle \z^{(-\infty,h+a)}, 1\rangle <\infty, \quad\N^{h}_0\hbox{ a.s.}
  \end{equation}
  and 
  \begin{equation}
  \label{key-tech2}
  \ve\,\inf\{k:Z^{\ve,k}=0\}\; \build{\la}_{\ve \to 0}^{\rm(d)} \;\inf\{a\geq 0:  \z^{(-\infty,h+a)} =0\}<\infty, \quad\N^{h}_0\hbox{ a.s.}
  \end{equation}
  We set 
  $$Z^{(\ve)}:=\sum_{k=0}^\infty Z^{\ve,2k}.$$

Next, from the special Markov property again and formula \eqref{law-max}, one obtains that, for every $k\geq 0$, the conditional
distribution of $n^k_\ve$ knowing $Z^{\ve,2k}$ is Poisson with parameter $\frac{3}{2\ve^2}\,Z^{\ve,2k}$. 
Simple Borel-Cantelli type arguments, using also \eqref{key-tech1} and \eqref{key-tech2}, now show that
$$\frac{n_\ve}{ \frac{3}{2\ve^2}\,Z^{(\ve)}} \build{\la}_{\ve\to 0}^{} 1$$
in $\N^{h}_0$-probability. So the proof of Lemma \ref{key-approx} will be
complete if we can verify that
\begin{equation}
\label{key-tech3}
2\ve\,Z^{(\ve)} \build{\la}_{\ve\to 0}^{} L^h,
\end{equation}
in $\N^{h}_0$-probability.

To this end, we note that, by \eqref{occu},
$$\frac{1}{\ve}\,\int_0^\sigma \mathrm{d}s\,\mathbf{1}_{[h-\ve,h]}(\wh W_s) \build{\la}_{\ve\to0}^{\rm a.s.} L^h$$
and we write
$$\int_0^\sigma \mathrm{d}s\,\mathbf{1}_{[h-\ve,h]}(\wh W_s)
= \int_0^\sigma \mathrm{d}s\,\mathbf{1}_{\{\tau_h(W_s)<\infty\}}\,\mathbf{1}_{[h-\ve,h]}(\wh W_s)
+ \int_0^\sigma \mathrm{d}s\,\mathbf{1}_{\{\tau_h(W_s)=\infty\}}\,\mathbf{1}_{[h-\ve,h]}(\wh W_s).$$
Recall that $B=(B_t)_{t\geq 0}$ stands for a linear Brownian motion starting from $0$ under the probability measure $P_0$,
and set $\theta_h=\inf\{t\geq 0: B_t=h\}$. By
the first-moment formula for the Brownian snake \cite[Proposition 4.2]{LZurich},
$$\N_0\Big(\int_0^\sigma \mathrm{d}s\,\mathbf{1}_{\{\tau_h(W_s)=\infty\}}\,\mathbf{1}_{[h-\ve,h]}(\wh W_s)\Big)
=E_0\Big[\int_0^{\theta_h} \mathrm{d}t\,\mathbf{1}_{[h-\ve,h]}(B_t)\Big] =\ve^2$$
where the last equality holds for $0<\ve\leq h$, by an application of a classical Ray-Knight theorem. By combining the previous observations, we see that we have also
\begin{equation}
\label{key-tech4}
\frac{1}{\ve}\,\int_0^\sigma \mathrm{d}s\,\mathbf{1}_{\{\tau_h(W_s)<\infty\}}\,\mathbf{1}_{[h-\ve,h]}(\wh W_s) \build{\la}_{\ve\to0}^{} L^h
\end{equation}
 in $\N^{h}_0$-probability. Then repeated applications of the special Markov property (in a way
 very similar to the proof of Lemma \ref{prelilemma}) show that we can write, $\N^h_0$-a.s.,
\begin{equation}
\label{key-tech5}
\int_0^\sigma \mathrm{d}s\,\mathbf{1}_{\{\tau_h(W_s)<\infty\}}\,\mathbf{1}_{[h-\ve,h]}(\wh W_s) 
 = U_\ve(Z^{(\ve)}),
 \end{equation}
 where $U_\ve= (U_\ve(r))_{r\geq 0}$ is a subordinator with no drift, whose L\'evy measure $\pi_\ve$
 is the ``law'' under $\N_0$ of 
 $$\int_0^\sigma \mathrm{d}s\,\mathbf{1}_{\{\tau_\ve(W_s)=\infty\}}\,\mathbf{1}_{[-\ve,0]}(\wh W_s).$$
 By simple scaling arguments, we have
 $$U_\ve(r)= \ve^4\,\wt U_\ve(\frac{r}{\ve^2})$$
 where $\wt U_\ve$ is a subordinator with L\'evy measure $\pi_1$, whose law does not depend on $\ve$. 
Using again the first-moment formula for the Brownian snake, we get that 
 $$\int y\,\pi_1(\mathrm{d}y)= \N_0\Big(\int_0^\sigma \mathrm{d}s\,\mathbf{1}_{\{\tau_1(W_s)=\infty\}}\,\mathbf{1}_{[-1,0]}(\wh W_s)\Big)
 = E_0\Big[\int_0^{\theta_1} \mathrm{d}t\, \mathbf{1}_{[-1,0]}(B_t)\Big] = 2.$$
 The law of large numbers then implies that
 \begin{equation}
 \label{key-tech6}
 \frac{U_\ve(Z^{(\ve)})}{\ve^2\,Z^{(\ve)}} = \frac{\ve^2}{Z^{(\ve)}}\,\wt U_\ve(\frac{Z^{(\ve)}}{\ve^2}) 
 \build{\la}_{\ve\to 0}^{} 2,
 \end{equation}
 in $\N^{h}_0$-probability. Our claim \eqref{key-tech3} now follows by combining
 \eqref{key-tech4},  \eqref{key-tech5} and  \eqref{key-tech6}. This completes the
 proof of Lemma \ref{key-approx}. \hfill$\square$
 
 \medskip
We finally explain how Theorem \ref{SBM} is derived from \cref{appro-LT} via the Brownian
snake construction of (historical) super-Brownian motion \cite[Chapter IV]{LZurich}, which we briefly recall below.

\medskip
\noi{\it Proof of Theorem \ref{SBM}.} We now argue
under the probability measure $\P_0$. For every $t\geq 0$, let $(\ell^t_s)_{s\geq 0}$ denote the
local time process at level $t$ of the reflected Brownian motion $(\zeta_s)_{s\geq 0}$. Fix $a>0$, and recall our
notation $\eta_a=\inf\{s\geq 0: \ell^0_s>a\}$. A historical super-Brownian motion $\mathbf{Y}$ starting from 
$a\,\delta_0$ can be obtained under $\P_0$ by setting, for every $t\geq 0$,
and every nonnegative measurable function $\Phi$ on $C([0,t],\R)$,
\begin{equation}
\label{SBMsnake}
\langle \mathbf{Y}_t,\Phi\rangle = \int_0^{\eta_a} d\ell^t_s\,\Phi(W_s),
\end{equation}
where the notation $d\ell^t_s$ refers to integration with respect to the increasing
function $s\la \ell^t_s$. In particular, if $\mathrm{supp}(\mathbf{Y}_t)$ stands for the 
topological support of $\mathbf{Y}_t$, we have a.s.~for every $t\geq 0$,
\begin{equation}
\label{support}
\mathrm{supp}(\mathbf{Y}_t) \subset \{W_s: s\in [0,\eta_a]\hbox{ and } \zeta_s=t\}.
\end{equation}
Conversely, formula \eqref{SBMsnake} implies that, a.s., for every $t\geq 0$, for every 
$s\in[0,\eta_a]$ such that $\zeta_s=t$ and $\sup\{\zeta_u:(s-\ve)_+\leq u\leq s+\ve\} > t$ for every $\ve>0$, one has
$W_s\in \mathrm{supp}(\mathbf{Y}_t)$ (in particular \eqref{support} is an equality if
$t$ is not a 
local maximum of $s\to\zeta_s$). Also note that, if $\mathbf{X}$ is the super-Brownian motion associated 
with $\mathbf{Y}$, the random measure $\int_0^\infty \mathrm{d}t\,\mathbf{X}_t$ coincides with the occupation
measure of $\wh W$ over the interval $[0,\eta_a]$. 

Write $N_{h,\ve}(a)$ for the number of upcrossing times of $W$
from $h$ to $h+\ve$ before time $\eta_a$. Before time $\eta_a$, there is only a finite number of excursions of
$W$ away from $0$ that hit level $h$, and obviously $N_{h,\ve}(a)$ is the sum of the upcrossing numbers
corresponding to each of these excursions. We can then apply \cref{appro-LT} to see
that $\ve^3 N_{h,\ve}(a)$ converges in probability to ($\frac{c_1}{2}$ times) the density at $h$ of the
occupation measure of $\wh W$ over $[0,\eta_a]$, which coincides with the local time of $\mathbf{X}$
at level $h$. 

From the previous considerations, the proof of Theorem \ref{SBM} will be complete if we
can verify, with the notation of this theorem, that $\mathscr{N}_{h,\ve}=N_{h,\ve}(a)$, a.s.~for every
fixed $h$ and $\ve$. In other words, we need to prove that upcrossing times 
of $\mathbf{Y}$ from $h$ to $h+\ve$ are in one-to-one correspondence with 
upcrossing times of $W$ from $h$ to $h+\ve$ before time $\eta_a$. 

Consider
an upcrossing time $r$
of $\mathbf{Y}$ from $h$ to $h+\ve$. By the definition, there exists $t>r$ and $w\in \mathrm{supp}(\mathbf{Y}_t)$
such that $r$ is an upcrossing time of the function $w$. From \eqref{support}, there exists $s\in[0,\eta_a]$
such that $\zeta_s=t$ and $W_s=w$. Set $\tilde s:=\sup\{u<s: \zeta_u=r\}$, so that
in particular $\zeta_{\tilde s}=r$. Then by the properties of
the Brownian snake the path $W_{\tilde s}$ coincides with the path $W_s=w$ restricted to $[0,r]$. Furthermore,
we have $\zeta_u>r$ for every $u\in(\tilde s,s]$ by construction, and it easily follows that 
$\tilde s$ is an upcrossing time of $W$ (if $\tilde r=\inf\{r'>r : w(r')=h+\ve\}$, take
$\tilde s'=\sup\{u<s: \zeta_u=\tilde r\}$, and note that $W_{\tilde s'}$ coincides with the
restriction of $w$ to $[0,\tilde r]$, so that the pair $(\tilde s,\tilde s')$ satisfies the properties
of the definition of an upcrossing time of the Brownian snake). 

By the previous discussion, for any  upcrossing time $r$
of $\mathbf{Y}$ from $h$ to $h+\ve$, we can construct an upcrossing time $\tilde s\in[0,\eta_a]$ of $W$ from
$h$ to $h+\ve$ such that $\zeta_{\tilde s}=r$. In fact, $\tilde s$ is uniquely determined by $r$: The point 
is that the quantities $\zeta_s$ when $s$ varies among upcrossing times of $W$ from $h$ to $h+\ve$
are distinct a.s. (recall that $h$ and $\ve$ are fixed). The latter property essentially follows from the fact that, if $B$ and $B'$
are two independent linear Brownian motions, the set of all left ends of excursion intervals of $B$ away from $h$
and the similar set for $B'$ are disjoint a.s. We omit some details here.

Clearly the mapping $r\to \tilde s$ is
one-to-one. It remains to verify that this mapping is also onto, and to this end it will suffice to
check that, for any upcrossing time $s$ of $W$ from
$h$ to $h+\ve$ before time $\eta_a$, $\zeta_s$ is an upcrossing time 
of $\mathbf{Y}$ from $h$ to $h+\ve$. Let $s\in[0,\eta_a]$ be an 
upcrossing time of $W$ from $h$ to $h+\ve$, and let $s'$ be the associated time. We already noticed that $\zeta_{s}$
is an upcrossing time of the function $W_{s'}$. If $t=\zeta_{s'}$, we have then $W_{s'}\in \mathrm{supp}(\mathbf{Y}_t)$,
by the observations following \eqref{support} and the fact that the time $s'$ cannot be a time of local maximum of $\zeta$ (this fact is a consequence of the strong Markov property of the Brownian snake, using the remark following Definition \ref{defup}).
Hence, we get that $\zeta_s$ is an upcrossing time of $\mathbf{Y}$ as desired. 

Finally, the mapping $r\to \tilde s$ is a bijection, and it follows that $\mathscr{N}_{h,\ve}=N_{h,\ve}(a)$, a.s. This completes the proof of
Theorem \ref{SBM}. \hfill $\square$

\smallskip
\rem Theorem \ref{SBM} can be extended to more general initial values of $\mathbf{X}$. In particular, the preceding proof shows that the result
still holds if $\mathbf{X}_0$ is supported on a compact interval $I$ and $h\notin I$. The convergence \eqref{intro-2} presumably holds for any
initial value $\mathbf{X}_0$ and any $h\in\R$. Proving this would however require some additional estimates.

 \section{Conditioned excursion measures}
 \label{sec-condi}
 
 In view of our applications to the Brownian map, we will now establish 
 certain conditional versions of \cref{appro-LT}. We first consider
 the probability measure $\N_0^{(1)}$ defined by
 $$\N_0^{(1)}=\N_0(\cdot \mid \sigma =1).$$
 Under $\N_0^{(1)}$, the lifetime process $(\zeta_s)_{s\geq 0}$ is a normalized
 Brownian excursion, and the conditional distribution of $(W_s)_{s\geq 0}$
 knowing $(\zeta_s)_{s\geq 0}$ remains the same as under $\N_0$. The definition 
 of upcrossing times of $W$ still makes sense under $\N_0^{(1)}$, and the local times
 $(L^h)_{h\in\R}$ are again well defined thanks to Theorem 2.1 in \cite{BMJ}.
 
 \begin{proposition}
 \label{appro-norma} 
 Let $h\in\R\backslash\{0\}$, and, for every $\ve>0$, let $N_{h,\ve}$ be the number of
 upcrossing times from $h$ to $h+\ve$. Then,
 $$\ve^3\, N_{h,\ve} \build{\la}_{\ve\to 0}^{} \frac{c_1}{2}\, L^h$$
 in $\N_0^{(1)}$-probability.
 \end{proposition}
 
 \rem It is very plausible that this result also holds for $h=0$, but we will leave this extension as an exercise
 for the reader, since it is not needed in our application to the Brownian map.
 
 \medskip
 \proof We rely on an absolute continuity argument to derive Proposition \ref{appro-norma}
 from \cref{appro-LT}. We fix $\eta>0$ and, on the event $\{\zeta_{1/2} > \eta\}$, we set
 \begin{align*}
 S_\eta&:=\inf\{s>\frac{1}{2} : \zeta_s=\eta\},\\
 R_\eta&:=\sup\{s<\frac{1}{2}:\zeta_s=\eta\}.
 \end{align*}
 If $\zeta_{1/2}\leq \eta$, we take $R_\eta=S_\eta=\frac{1}{2}$.
 From standard facts about Brownian excursions (we omit a few details here), one 
 easily checks that the law of the process
 $$(\zeta_{(R_\eta+s)\wedge S_\eta} - \eta)_{s\geq 0}$$
 under the conditional probability measure $\N_0^{(1)}(\cdot\mid \zeta_{1/2} > \eta)$
 is absolutely continuous with respect to the It\^o measure of Brownian excursions. 
 Note that, for $R_\eta\leq s\leq S_\eta$, we have $W_s(\eta)=W_{R_\eta}(\eta)=\wh W_{R_\eta}$
 by the properties of the Brownian snake. Still on the event $\{\zeta_{1/2} > \eta\}$, we define a path-valued process
 $W^\eta=(W^\eta_s)_{s\geq 0}$ by setting, for every $s\geq 0$,
 $$W^\eta_s := W_{(R_\eta+s)\wedge S_\eta} (\eta +t) - \wh W_{R_\eta}\;,\quad 0\leq t \leq 
\zeta^\eta_s:= \zeta_{(R_\eta+s)\wedge S_\eta} -\eta.$$
Then the law of $(W^\eta_s)_{s\geq 0}$ under $\N_0^{(1)}(\cdot\mid \zeta_{1/2} > \eta)$ is absolutely continuous with respect to $\N_0$. Furthermore,
the process $(W^\eta_s)_{s\geq 0}$ is independent of $H_\eta:=\wh W_{R_\eta}$ under
the same probability measure. These facts are simple consequences 
of the properties of the Brownian snake.

Let $h>0$ (the case $h<0$ is treated in a similar way). We can choose $\eta>0$ small enough in such a way that,
$$\N_0^{(1)}(\zeta_{1/2} \leq \eta) + \N_0^{(1)}\Big( \zeta_{1/2} >\eta; \Big(\sup_{s\leq R_\eta} \wh W_s\Big)
\vee \Big(\sup_{s\geq S_\eta} \wh W_s\Big) \geq  h\Big)$$
is arbitrarily small. On the other hand, on the event
$$A_\eta:= \{\zeta_{1/2} > \eta\} \cap \Big\{ \Big(\sup_{s\leq R_\eta} \wh W_s\Big)
\vee \Big(\sup_{s\geq S_\eta} \wh W_s\Big) < h\Big\}$$
simple considerations give, with an obvious notation,
$$L^h = L^{h-H_\eta}(W^\eta)$$
and, for every $\ve>0$,
$$N_{h,\ve} = N_{h-H_\eta,\ve}(W^\eta).$$
Because the law of $W^\eta$ is absolutely continuous with respect
to $\N_0$, and using also the fact that $H_\eta$ is independent of $W_\eta$ under $\N_0^{(1)}(\cdot\mid \zeta_{1/2} > \eta)$,
we can use \cref{appro-LT} to obtain that the convergence
$$\ve^3\, N_{h-H_\eta,\ve}(W^\eta) \build{\la}_{\ve\to0}^{} \frac{c_1}{2}\,L^{h-H_\eta}(W^\eta),$$
holds in probability under $\N_0^{(1)}(\cdot\mid \zeta_{1/2} > \eta)$. The result of Proposition 
\ref{appro-norma} follows from the preceding considerations.  \endproof

We finally give the analog of Proposition 
\ref{appro-norma} for the Brownian snake ``conditioned to stay positive''. It is proved in 
\cite{LGW} that the conditional measures $\N^{(1)}_0(\cdot \mid \inf\{\wh W_s:s\geq 0\} >-\delta)$
converge as $\delta\downarrow 0$ to a limit, which is denoted by $\overline\N^{(1)}_0$. This limiting measure
can also be constructed directly as the law under $\N^{(1)}_0$ of the Brownian snake ``re-rooted'' at its
minimum. Let us describe this construction
 (see \cite{LGW} for more details). We argue under the measure $\N^{(1)}_0$. Fix
 $r\in[0,1]$, and  set, for every $s\in[0,1]$,
 $$\zeta^{[r]}_s:=d_\zeta(p_\zeta(r),p_\zeta(r\oplus s))=\zeta_r +\zeta_{r\oplus s} - 2\min_{u\in [r\wedge (r\oplus s),r\vee (r\oplus s)]} \zeta_u$$
 with the notation $r\oplus s= r+s$ if $r+s\leq 1$, and $r\oplus s=r+s-1$ if $r+s>1$. Also set
 $\zeta^{[r]}_s=0$ if $s>1$. Then, the tree $\t_{\zeta^{[r]}}$ is identified isometrically with the
 tree $\t_\zeta$ re-rooted at $p_\zeta(r)$, via the mapping $p_{\zeta^{[r]}}(s)\la p_\zeta(r\oplus s)$. We also introduce a path-valued process
 $W^{[r]}$ such that the associated lifetime process is $\zeta^{[r]}$: We first set
 $\wh W^{[r]}_s:= \wh W_{r\oplus s} - \wh W_r$, for every $s\in[0,1]$, and we then define
 the path $W^{[r]}_s$ by saying that, for every $t\in[0,\zeta^{[r]}_s]$, $ W^{[r]}_s(t) = \wh W^{[r]}_u$
 if $u\in[0,1]$ is such that $p_{\zeta^{[r]}}(u)$ is the (unique) ancestor of $p_{\zeta^{[r]}}(s)$
 at distance $t$ from the root in the tree $\t_{\zeta^{[r]}}$. The invariance of the Brownian snake under
 re-rooting (see formula (3) in \cite{LGW}) asserts that $(W^{[r]}_s)_{0\leq s\leq 1}$ has the same 
 distribution as $(W_s)_{0\leq s\leq 1}$ under $\N^{(1)}_0$.

 Of course the preceding invariance property fails if we allow $r$ to be random. We let $s_*$ be the almost surely
 unique element of $[0,1]$ such that
 $$\wh W_{s_*}= \min\{\wh W_s: s\in[0,1]\}.$$
By \cite[Theorem 1.2]{LGW}, the process $W^{[s_*]}$ is distributed according to
 $\overline\N^{(1)}_0$. Furthermore, $s_*$ is uniformly distributed over $[0,1]$, and 
 $s_*$ and $W^{[s_*]}$ are independent under $\N^{(1)}_0$. The latter two properties are
 straightforward consequences of the invariance under (deterministic) re-rooting. 
 
 We write $W_*=\wh W_{s_*}$ to simplify notation. In a way similar to the discussion
 before Lemma \ref{fini-upcro}, we assign the spatial location $\Gamma^{[s_*]}_v=\wh W^{[s*]}_s$ to the vertex $v=
 p_{\zeta^{[s*]}}(s)$ of
 $\t_{\zeta^{[s_*]}}$, for every $s\in[0,1]$, and,
 modulo the identification of $\t_{\zeta^{[s_*]}}$ 
with $\t_\zeta$ re-rooted at $p_\zeta(s_*)$, we have $\Gamma^{[s_*]}_v=\Gamma_v - W_{*}$
for every $v$. 
 
Furthermore, it is clear that the definition of the local times $L^h$ still makes sense
under $\overline\N^{(1)}_0$: Just note that the occupation measure of $\wh W^{[s*]}$
concides with the occupation measure of $\wh W$
shifted by $-W_{*}$. 

\begin{proposition}
\label{appro-condi}
Let $h>0$. The convergence of Proposition \ref{appro-norma} also holds in $\overline\N^{(1)}_0$-probability.
\end{proposition}

\proof 
We fix $\kappa>0$. We can choose $\alpha\in(0,1/4)$ such that 
\begin{equation}
\label{conditech1}
\overline\N^{(1)}_0\Big(\sup_{s\in [0,2\alpha]\cup [1-2\alpha,1]} |\wh W_s|> \frac{h}{4}\Big) < \kappa.
\end{equation}
Then, we choose $\eta>0$ such that
\begin{equation}
\label{conditech2}
\overline\N^{(1)}_0\Big(\inf_{s\in [\frac{\alpha}{2},1-\frac{\alpha}{2}]} \zeta_s \leq 2\eta\Big) < \kappa.
\end{equation}
Finally, recalling that $(\zeta_s)_{0\leq s\leq 1}$ is distributed under $\N^{(1)}_0$ as a normalized
Brownian excursion, we can choose $\delta\in (0,\frac{\alpha}{2})$ such that
$$\N^{(1)}_0\Big(\sup_{s\leq \delta} \zeta_s >\eta\Big) < \kappa\,\delta.
$$
Since $s_*$ is uniformly distributed over $[0,1]$, this last bound also implies
\begin{equation}
\label{conditech3}
\N^{(1)}_0\Big(\sup_{s\leq \delta} \zeta_s >\eta\,\Big|\, s_*<\delta\Big) < \kappa.
\end{equation}
From the results recalled before the statement of the proposition, and in particular
the fact that $s_*$ and $W^{[s_*]}$ are independent under $\N^{(1)}_0$, we obtain that,
under the conditional probability measure $\N^{(1)}_0(\cdot|\, s_*<\delta)$, the process
$W^{[s_*]}$ is distributed according to 
$\overline\N^{(1)}_0$, so that we can apply the bounds \eqref{conditech1} and \eqref{conditech2}
to this process. Combining the bounds \eqref{conditech1}, \eqref{conditech2} and \eqref{conditech3},
we see that, except on a set of $\N^{(1)}_0(\cdot|\, s_*<\delta)$-measure smaller than $3\kappa$, we have
\begin{enumerate}
\item[(i)] $\forall s\in [\frac{\alpha}{2}, 1-\frac{\alpha}{2}],\; \zeta^{[s_*]}_s > 2\eta$;
\item[(ii)] $\forall s\in [0,2\alpha]\cup [1-2\alpha,1],\; |\wh W^{[s_*]}_s| \leq \frac{h}{4}$;
\item[(iii)] ${\displaystyle \sup_{s\leq \delta} \zeta_s \leq \eta}$.
\end{enumerate}
Now recall the definition of $\zeta^{[s_*]}$ and $\wh W^{[s_*]}$ in terms of the pair 
$(\zeta,\wh W)$. Using also the fact that $\delta<\frac{\alpha}{2}$, we see that (i)--(iii)
imply, on the event $\{s_*<\delta\}$,
\begin{enumerate}
\item[(i)'] $\forall s\in [\alpha,1-\alpha],\;\zeta_s > \eta$;
\item[(ii)'] $\forall s\in[0,\alpha]\cup[1-\alpha,1],\; |\wh W_s| \leq \frac{h}{4}$.
\end{enumerate}
Recall the notation $R_\eta,S_\eta,W^\eta,H_\eta$ introduced in the proof of Proposition \ref{appro-norma}. 
Obviously (i)' implies that $R_\eta<\alpha$ and $S_\eta>1-\alpha$. Therefore we can summarize 
the preceding discussion by saying that, except on a set of $\N^{(1)}_0(\cdot|\, s_*<\delta)$-measure smaller than $3\kappa$,
we have both the properties (i)' and (iii) above, and
\begin{equation}
\label{conditech4}
\sup_{s\in[0,R_\eta]\cup [S_\eta,1]} |\wh W_{s}| \leq \frac{h}{4}.
\end{equation}
Write $A_{\delta,\eta}$ for the intersection of the event $\{s_*<\delta\}$ with the event where (i)', (iii) and \eqref{conditech4} hold, and use the obvious notation
$N_{h,\ve}(W^{[s_*]})$ and $L^h(W^{[s_*]})$ for respectively the upcrossing numbers and the local times of
$W^{[s_*]}$. Also note that (iii) forces $\delta \leq R_\eta$. Then,
\begin{align}
\label{bigdisplay}
&\overline\N^{(1)}_0\Big(|\ve^3N_{h,\ve} - \frac{c_1}{2} L^h| \wedge 1\Big)\nonumber\\
&\quad= \N^{(1)}_0 \Big( |\ve^3N_{h,\ve}(W^{[s_*]})  - \frac{c_1}{2} L^h(W^{[s_*]})| \wedge 1\,\Big|\,s_*<\delta\Big)\nonumber\\
&\quad\leq \N^{(1)}_0 \Big(\mathbf{1}_{A_{\delta,\eta}}\,|\ve^3N_{h,\ve}(W^{[s_*]}) - \frac{c_1}{2} L^h(W^{[s_*]})| \wedge 1\,\Big|\,s_*<\delta\Big) + 3\kappa\nonumber\\
&\quad= \N^{(1)}_0 \Big(\mathbf{1}_{A_{\delta,\eta}}\,|\ve^3N_{h+W_*,\ve} - \frac{c_1}{2} L^{h+W_*}| \wedge 1\,\Big|\,s_*<\delta\Big) + 3\kappa\nonumber\\
&\quad= \N^{(1)}_0 \Big(\mathbf{1}_{A_{\delta,\eta}}\,|\ve^3N_{h+W_*-H_\eta,\ve} (W^\eta)- \frac{c_1}{2} L^{h+W_*-H_\eta}
(W^\eta)| \wedge 1\,\Big|\,s_*<\delta\Big) + 3\kappa.
\end{align}
In the fourth line of the preceding display, we use \eqref{conditech4} (and the fact that $s_*<\delta\leq R_\eta$)
to verify that $N_{h,\ve}(W^{[s_*]})= N_{h+W_*,\ve}$ and $L^h(W^{[s_*]})= L^{h+W_*}$ on the event
$A_{\delta,\eta}$. In particular, the simplest way to obtain the identity $N_{h,\ve}(W^{[s_*]})= N_{h+W_*,\ve}$ is
to use the interpretation of upcrossing times in terms of vertices of the tree $\t_\zeta$ (see the discussion before Lemma \ref{fini-upcro}), observing that $\t_{\zeta^{[s_*]}}$ 
is identified with $\t_\zeta$ re-rooted at $p_\zeta(s_*)$ and that, modulo this identification, a vertex
$v$ of $\t_\zeta$ such that $\Gamma_v=h$ has, on the event $A_{\delta,\eta}$, the same descendants in $\t_\zeta$
and in $\t_{\zeta^{[s_*]}}$. Similarly, in the last equality of \eqref{bigdisplay}, we use \eqref{conditech4} to replace $N_{h+W_*,\ve}$ and $L^{h+W_*}$
by $N_{h+W_*-H_\eta,\ve} (W^\eta)$ and $L^{h+W_*-H_\eta}
(W^\eta)$ respectively. 

Clearly, in the last line of \eqref{bigdisplay}, we can replace $W_*$ by
$$W_*^{(R_\eta)}:= \min\{\wh W_s:0\leq s\leq R_\eta\}.$$
Under the probability measure $\N^{(1)}_0$, if we condition on $R_\eta$ and $S_\eta$, $W^\eta$ becomes {\it independent}
of the pair $(W_*^{(R_\eta)},H_\eta)$, and is distributed as a Brownian snake excursion with duration $S_\eta-R_\eta$.
Therefore we can apply Proposition \ref{appro-norma} to see that
$$\lim_{\ve\to 0} \N^{(1)}_0 \Big(\mathbf{1}_{\{R_\eta\leq \alpha\leq 1-\alpha\leq S_\eta,
|W_*^{(R_\eta)}-H_\eta|\leq \frac{h}{2}\}}\,\Big|\ve^3N_{h+W_*^{(R_\eta)}-H_\eta,\ve} (W^\eta)- \frac{c_1}{2} L^{h+W_*^{(R_\eta)}-H_\eta}
(W^\eta)\Big| \wedge 1\Big) =0.$$
Noting that $A_{\delta,\eta}\subset \{R_\eta\leq \alpha\leq 1-\alpha\leq S_\eta,
|W_*^{(R_\eta)}-H_\eta|\leq \frac{h}{2}\}$, we now conclude from \eqref{bigdisplay} that
$$\limsup_{\ve \to 0} \;\overline\N^{(1)}_0\Big(|\ve^3N_{h,\ve} - \frac{c_1}{2} L^h| \wedge 1\Big)\leq 3\kappa$$
and since $\kappa$ was arbitrary this completes the proof. \endproof

\section{Application to the Brownian map}
\label{sec-applimap}

Let us recall the construction of the Brownian map $(\bm_\infty, D)$ from the Brownian snake
$(W_s)_{0\leq s\leq 1}$. In the following presentation,
we argue under the probability measure $\ov\N^{(1)}_0$. The fact that $\ov\N^{(1)}_0$ 
coincides with the law under $\N^{(1)}_0$ of the Brownian snake ``re-rooted at its 
minimum'' (as explained above) shows that this presentation is equivalent to the one given 
in \cite{LUniqueness} or in \cite{LGM}  . We can define the collection $(\Gamma_a)_{a\in\t_\zeta}$
under the probability measure $\ov\N^{(1)}_0$ by setting $\Gamma_a=\wh W_s$ if $p_\zeta(s)=a$, exactly as we did  under $\N_0$
in Section \ref{sec-upcrossing}. Note that $\Gamma_{\rho_\zeta}=0$ and 
$\Gamma_a\geq 0$ for every $a\in\t_\zeta$, $\ov\N^{(1)}_0$ a.s. We now
interpret $\Gamma_a$ as a label assigned to the vertex $a$. For every $a,b\in\t_\zeta$, we then set
$$D^\circ(a,b):=\Gamma_a+\Gamma_b - 2\max_{
\begin{subarray}{c}
s,s'\in[0,1]\\
 p_\zeta(s)=a,p_\zeta(s')=b
 \end{subarray}}\Big( \min\{ \wh W_r : s\wedge s'\leq r \leq s\vee s'\}\Big)$$
 and 
 $$D(a,b):= \inf\Big\{ \sum_{i=1}^p D^\circ(a_{i-1},a_i)\Big\}$$
 where the infimum is over all choices of the integer $p\geq 1$ and of the elements
 $a_0,a_1,\ldots,a_p$ of $\t_\zeta$ such that $a_0=a$ and $a_p=b$. Then
 $D$ is  a pseudo-metric on $\t_\zeta$, and we consider the associated 
 equivalence relation $\approx$: if $a,b\in\t_\zeta$,
 $$a\approx b\hbox{ \ if and only if \ }D(a,b)=0\,.$$
 One can prove that this property is also equivalent to $D^\circ(a,b)=0$.
 Informally, this means that $a$ and $b$ have the same label, and that one 
 can go from $a$ to $b$ moving ``around'' the tree and encountering only vertices 
 with a larger label. 
  
 The Brownian map is the quotient space $\bm_\infty:=\t_\zeta\,/\! \approx$, which is equipped
 with the metric induced by $D$, for which we keep the same notation $D$. We write $\Pi$
 for the canonical projection from $\t_\zeta$ onto $\bm_\infty$. The projection $\Pi$ is continuous
 (see \cite[Section 2.5]{LGeodesic}). We will use the following lower 
 bound \cite[Corollary 3.2]{LGeodesic}: For every $a,b\in\t_\zeta$,
 \begin{equation}
 \label{lowerbd}
 D(\Pi(a),\Pi(b))\geq \Gamma_a + \Gamma_b - 2\,\inf_{c\in \llbracket a,b \rrbracket} \Gamma_c.
 \end{equation}

The distinguished point of the Brownian map is $\rho=\Pi(\rho_\zeta)$, where we recall that $\rho_\zeta=p_\zeta(0)$ is the root of $\t_\zeta$.
 We have then $D(\rho,\Pi(a))= \Gamma_a$ for every $a\in \t_\zeta$. 
 The volume measure $\lambda$ on $\bm_\infty$ is the image of 
 Lebesgue measure on $[0,1]$ under $\Pi\circ p_\zeta$. 
 
 In the proof of Theorem \ref{connec-compo}, we will need the following lemma. We say
 that $h\in\R$ is a local minimum of $\w\in\W$ if there exists $t\in(0,\zeta_{(\w)})$ and $\beta>0$,
 with $(t-\beta,t+\beta)\subset (0,\zeta_{(\w)})$,
 such that $\w(t)=h$ and $\w(t')\geq h$ for every $t'\in(t-\beta,t+\beta)$. 
 
 \begin{lemma}
 \label{localmin}
 Let $h>0$. Then $\ov\N^{(1)}_0$ a.s. for every $s\in[0,1]$, $h$ is not a local minimum of $W_s$.
 \end{lemma}
 
\proof If we replace $\ov\N^{(1)}_0$ by $\N^{(1)}_0$ in the statement of the lemma, the proof is easy, by an argument
already explained at the end of
the proof of Lemma \ref{prelilemma2}. To get
 the precise statement of the lemma, we need to verify that $h+W_*$ is not a local minimum of one
 of the paths $W_s$,  $\N^{(1)}_0$ a.s.  The fact that $h+W_*$ is random, and of course not independent
 of the paths $W_s$, makes the proof a little harder. Still one can 
 use arguments very similar to the proof of Proposition \ref{appro-condi}, conditioning on
 the event $\{s_*<\delta\}$ and replacing $W_*$ by $W_*^{(R_\eta)}$ (where $\delta$ and
 $\eta$ are chosen as in the latter proof): Except on a set of small
 probability, one can then concentrate on the paths $W_s$ for $s\in[R_\eta,S_\eta]$, or more precisely on 
 the paths $W^\eta_s$ for $s\in[0,S_\eta-R_\eta]$,
 and use the same independence property as in the end of the proof of Proposition \ref{appro-condi} 
 to conclude. We leave the details to the reader. \endproof

 \smallskip
 \noi{\it Proof of Theorem \ref{connec-compo}.} It easily follows from the formula  $D(\rho,\Pi(a))= \Gamma_a$ and our definition
 of the volume measure $\lambda$ that the profile of distances $\Delta$ coincides with the occupation measure
 of the (conditioned) Brownian snake. Consequently, $\lambda$ has a continuous density 
 $(\mathbf{L}^h)_{h\geq 0}$ and $\mathbf{L}^h=L^h$, $\ov\N^{(1)}_0$ a.s. We then claim that,
 for every fixed $h>0$ and $\ve>0$,
 $$\mathbf{N}_{h,\ve}=N_{h,\ve},\quad 
 \ov\N^{(1)}_0\hbox{ a.s. }$$
Once the claim is proved, the statement of the theorem follows from Proposition \ref{appro-condi}. 

Say that $a\in\t_\zeta$ is an $(h,\ve)$-upcrossing vertex if $a=p_\zeta(s)$ where $s$ is an
upcrossing time of the Brownian snake from $h$ to $h+\ve$. This is equivalent to saying that
$\Gamma_a=h$ and $a$ has a descendant $b$ such that $\Gamma_b=h+\ve$
and $\Gamma_c>h$ for every $c\in \llbracket a,b\rrbracket\backslash\{a\}$. Note that
we have then $D(\rho,a)=\Gamma_a=h$. To prove our claim, we verify that
$(h,\ve)$-upcrossing vertices are in one-to-one correspondence with connected
components of 
 $B_h(\rho)^c$ that intersect $B_{h+\ve}(\rho)^c$.
 
 Let $a$ be an $(h,\ve)$-upcrossing vertex. We define $C_a$ as the set of all vertices $b\in\t_\zeta$
 such that $b$ is a descendant of $a$ and $\Gamma_c\geq h$ for every $c\in\llbracket a,b\rrbracket$. Note that if $b\in C_a$, then the whole segment $\llbracket a,b\rrbracket$ is contained
 in $C_a$. It follows that $C_a$ is (path-)connected, and it is also easy to check that $C_a$ is a closed
 subset of $\t_\zeta$.
 Furthermore the fact that $a$ is an $(h,\ve)$-upcrossing vertex ensures that $C_a$
 contains (at least) one vertex $a'$ such that $\Gamma_{a'}= a+\ve$. 
 
 To simplify notation, set $\t_\zeta^{\geq h}:=\{b\in\t_\zeta: \Gamma_b\geq h\}$. We next verify that
 $C_a$ is a connected component of $\t_\zeta^{\geq h}$. To this end, we set for every $\delta>0$,
 $$O_\delta:=\Big\{b\in  \t_\zeta^{\geq h}: \inf_{c\in\llbracket a,b \rrbracket} \Gamma_c > h-\delta\Big\}.$$
 It is easy to verify that $O_\delta$ is open in $\t_\zeta^{\geq h}$.  The set $O_\delta$ is also closed. In fact,
 if $(b_n)$ is a sequence in $O_\delta$ that converges to $b$, and if, for every $n$, $c_n$
 is the unique vertex of $\t_\delta$ such that $\llbracket a,b \rrbracket\cap \llbracket a,b_n \rrbracket= 
 \llbracket a,c_n \rrbracket$, then we must have $d_\zeta(c_n,b)\la 0$ as $n\to\infty$, and since
 $\Gamma_b=\lim \Gamma_{b_n} \geq h$, it follows that for $n$ large enough we have
 $$\inf_{c\in\llbracket a,b \rrbracket} \Gamma_c \geq\Big( \inf_{c\in\llbracket a,b_n \rrbracket} \Gamma_c \Big)
 \wedge (h-\frac{\delta}{2})> h-\delta,$$
 so that $b\in O_\delta$ as desired. We then observe that
\begin{equation}
\label{inter-pro}
C_a=\bigcap_{\delta>0} O_\delta.
\end{equation}
 Indeed, let $b\in C_a^c$. If $b$ is a descendant of $a$,  we must
 have 
 $$\inf_{c\in\llbracket a,b\rrbracket} \Gamma_c <h$$
yielding that $b\in O_\delta^c$ as soon as $\delta$ is small enough. If $b$ is not a descendant of $a$, we observe 
that $\llbracket a,b\rrbracket \cap \llbracket \rho,a\rrbracket=\llbracket \tilde a, a\rrbracket$, for some
$\tilde a\in \llbracket \rho,a\rrbracket$ such that $\tilde a\not = a$. Now recall that $a$ has a descendant $a'$
such that $\Gamma_{a'}=h+\ve$ and $\Gamma_c\geq h$ for every $h\in \llbracket a,a'\rrbracket$. From Lemma \ref{localmin},
the values of $\Gamma$ along $\llbracket \rho,a'\rrbracket$ cannot have a local minimum equal to $h$, and
we again obtain that
$$\inf_{c\in\llbracket a,b\rrbracket} \Gamma_c\leq \inf_{c\in\llbracket \tilde a,a\rrbracket} \Gamma_c <h.$$ 
Since we know that all sets $O_\delta$ are both open and closed in $\t_\zeta^{\geq h}$, \eqref{inter-pro}
implies that $C_a$ is a connected component of $\t_\zeta^{\geq h}$. 

Now let $\mathcal{C}_a:=\Pi(C_a)$. The preceding considerations  imply that
$\mathcal{C}_a$ is a connected component of $\Pi(\t_\zeta^{\geq h})=B_h(\rho)^c$. Let us explain this. By the continuity 
of $\Pi$, $\mathcal{C}_a$ is (path-)connected and closed in $\bm_\infty$. From \eqref{inter-pro}
and a simple compactness argument, we get that
\begin{equation}
\label{inter-pro2}
\mathcal{C}_a=\bigcap_{\delta>0} \Pi(O_\delta).
\end{equation}
The sets $\Pi(O_\delta)$ are closed by the continuity of $\Pi$. Let us prove that they are also open 
in  $B_h(\rho)^c$. Since we already know that $\t_\zeta^{\geq h}\backslash O_\delta$
is closed, this will  follow from the equality
$$B_h(\rho)^c\,\backslash\, \Pi(O_\delta)= \Pi(\t_\zeta^{\geq h}\backslash O_\delta).$$
In this equality, the inclusion $\subset$ is obvious. To prove the reverse inclusion, 
we need to verify that $ \Pi(O_\delta)\cap \Pi(\t_\zeta^{\geq h}\backslash O_\delta)=\varnothing$. Let
$b\in O_\delta$ and $\tilde b\in \t_\zeta^{\geq h}\backslash O_\delta$.
We have then $\Gamma_b\geq h$, $\Gamma_{\tilde b}\geq h$, and  
$$\inf_{c\in\llbracket a,b \rrbracket} \Gamma_c > h-\delta\;,\quad
 \inf_{c\in\llbracket a,\tilde b \rrbracket} \Gamma_c \leq h-\delta.$$
 Since $\llbracket b,\tilde b \rrbracket\supset \llbracket a,\tilde b \rrbracket\backslash \llbracket a,b \rrbracket$, it follows that
 $$ \inf_{c\in\llbracket b,\tilde b \rrbracket} \Gamma_c\leq h-\delta< \Gamma_b$$
 so that $\Pi(b)\not = \Pi(\tilde b)$ by \eqref{lowerbd}. We have thus proved 
 that the sets $\Pi(O_\delta)$ are both closed and open in $B_h(\rho)^c$, and \eqref{inter-pro2} now implies that
 $\mathcal{C}_a$ is a connected component of $B_h(\rho)^c$. 
 
 Summarizing, with each $(h,\ve)$-upcrossing vertex $a$ we can associate a connected component $\mathcal{C}_a$
 of $B_h(\rho)^c$
 that intersects $B_{h+\ve}(\rho)^c$. If $a$ and $ \tilde a$ are two distinct  $(h,\ve)$-upcrossing vertices, 
 we have $\mathcal{C}_a \not =\mathcal{C}_{\tilde a}$ because otherwise $\tilde a$ would be a descendant of $a$
 and $a$ would be a descendant of $\tilde a$, which is only possible if $a=\tilde a$. So it only remains to show that any connected component
  of $B_h(\rho)^c$ that intersects $B_{h+\ve}(\rho)^c$ is of this form. Let $\mathcal{C}$ be such a 
  connected component and let $x\in\mathcal{C}\cap B_{h+\ve}(\rho)^c$.
  Choose $b\in \t_\zeta$ such that $\Pi(b)=x$. Then $\Gamma_b\geq h+\ve$. By a continuity
  argument, there exists a unique vertex $a\in \llbracket \rho,b\rrbracket$ such that $\Gamma_a=h$
  and $\Gamma_c>h$ for every $c\in \llbracket a,b\rrbracket\backslash\{a\}$. Then $a$ is
  an $(h,\ve)$-upcrossing vertex, and $\mathcal{C}= \mathcal{C}_a$. This completes the proof.
  \hfill $\square$
  
  \medskip
  \rem With the notation of the preceding proof, set
  $C^\circ_a:=\{b\in C_a: \Gamma_b>h\}$ and $\mathcal{C}^\circ_a:=\Pi(C^\circ_a)$,
  for every $(h,\ve)$-upcrossing vertex $a$. Then the sets
  $\mathcal{C}^\circ_a$ are open in $\bm_\infty$ and these sets, when $a$ varies among all $(h,\ve)$-upcrossing vertices, are exactly those
  connected components of the complement of the closed ball $\ov B_h(\rho)$ that intersect $B_{h+\ve}(\rho)^c$. 
  So we could have stated Theorem \ref{connec-compo} in terms of connected components of the
  open set $\ov B_h(\rho)^c$, and the preceding proof would have been a little simpler. We chose to deal
  with connected components of the complement of the open ball mainly in view of the connection with 
  the Brownian cactus \cite[Section 2.5]{CLM}. As a final remark, it is not hard to verify that the boundary of $\mathcal{C}^\circ_a$,
  $\partial \mathcal{C}^\circ_a =\Pi(\{b\in C_a:\Gamma_b=h\})$, is a simple loop in $\bm_\infty$. By Jordan's theorem,
  all sets $\mathcal{C}^\circ_a$ are homeomorphic to the open unit disk. 
  
  \bigskip
  
  \noi{\bf Acknowledgement.} I thank Nicolas Curien and Gr\'egory Miermont for useful conversations about this
  work, which is a continuation of our previous article in collaboration \cite{CLM}.

\end{document}